\documentclass[twoside]{article}

\usepackage[accepted]{aistats2025}
\usepackage{url}            
\usepackage{booktabs}       
\usepackage{amsfonts}       
\usepackage{nicefrac}       
\usepackage{microtype}      
\usepackage{xcolor}         
\usepackage{algorithm}
\usepackage[]{algpseudocode}
\usepackage{multirow}
\usepackage{graphicx}
\usepackage{subfigure}
\usepackage{paralist}
\usepackage[font=small]{caption}

\usepackage{amssymb,amsmath,amsthm,dsfont}
\usepackage{bm}

\providecommand{\lin}[1]{\ensuremath{\left\langle #1 \right\rangle}}

\providecommand{\norm}[1]{\left\lVert#1\right\rVert_2}

\providecommand{\refLE}[1]{\ensuremath{\stackrel{(\ref{#1})}{\leq}}}
\providecommand{\refEQ}[1]{\ensuremath{\stackrel{(\ref{#1})}{=}}}

\providecommand{\R}{\mathbb{R}} 

\DeclareMathOperator{\E}{{\mathbb E}}
\providecommand{\E}[1]{{\mathbb E}\left.#1\right. }     
\providecommand{\Eb}[1]{\E \left[#1\right] }             
\providecommand{\EE}[2]{\E_{#1} \! #2 }      
\providecommand{\EEb}[2]{\E_{#1}\!\! \left[#2\right] } 

\DeclareMathOperator*{\argmin}{arg\,min}

\providecommand{\0}{\mathbf{0}}

\renewcommand{\gg}{\mathbf{g}}

\let\lll\ll
\renewcommand{\ll}{\mathbf{l}}

\providecommand{\xx}{\mathbf{x}}
\providecommand{\yy}{\mathbf{y}}

\providecommand{\mA}{\mathbf{A}}

\providecommand{\mH}{\mathbf{H}}

\providecommand{\mQ}{\mathbf{Q}}

\providecommand{\cL}{\mathcal{L}}
\providecommand{\cM}{\mathcal{M}}

\providecommand{\cO}{\mathcal{O}}

\RequirePackage[colorinlistoftodos,bordercolor=orange,backgroundcolor=orange!20,linecolor=orange,textsize=scriptsize]{todonotes}
\providecommand{\mycomment}[3]{\todo[inline,caption={},color=#3!20]{\textbf{#1: }#2}}%
\providecommand{\inlinecomment}[3]{%
{\color{#1}#2: #3}}%
\newcommand\commenter[2]%
{%
\expandafter\newcommand\csname i#1\endcsname[1]{\inlinecomment{#2}{#1}{##1}} 
\expandafter\newcommand\csname #1\endcsname[1]{\mycomment{#1}{##1}{#2}} 
}

\definecolor{mydarkblue}{rgb}{0,0.08,0.45}
\usepackage{hyperref} 
\usepackage[capitalize,noabbrev]{cleveref}

\newtheorem{lemma}{Lemma}

\newtheorem{assumption}{Assumption}
\crefname{assumption}{Assumption}{Assumptions}
\newtheorem{theorem}{Theorem}

\newtheorem{strongerassumption}{Assumption}
\crefalias{strongerassumption}{assumption}

\usepackage{thmtools}
\declaretheoremstyle[
spaceabove=\topsep, spacebelow=\topsep,
headfont=\normalfont\scshape, 
notefont=\mdseries, notebraces={(}{)},
bodyfont=\normalfont,
postheadspace=1em,
qed=\qedsymbol
]{mythmstyle}
\declaretheorem[style=mythmstyle]{remark}

\usepackage{url}

\usepackage{breakcites}
\usepackage{nicefrac}
\usepackage[normalem]{ulem}

\usepackage{enumitem}
\setlist[itemize]{leftmargin=20pt, itemsep=1pt, topsep=2pt}
\setlist[enumerate]{leftmargin=24pt, itemsep=1pt, topsep=2pt}

\usepackage{hyperref} 
\hypersetup{
hypertexnames=false, 
linkcolor=blue,
citecolor=cyan,
filecolor=green
}

\providecommand{\T}{^{\textnormal{T}}} 

\renewcommand{\epsilon}{\varepsilon}

\renewcommand{\phi}{\varphi}

\renewcommand{\Gamma}{\varGamma}

\renewcommand{\Delta}{\varDelta}

\renewcommand{\Theta}{\varTheta}
\let\oldLambda\Lambda
\renewcommand{\Lambda}{\varLambda}

\renewcommand{\Xi}{\varXi}

\renewcommand{\Pi}{\varPi}

\renewcommand{\Sigma}{\varSigma}

\renewcommand{\Upsilon}{\varUpsilon}

\renewcommand{\Phi}{\varPhi}

\renewcommand{\Psi}{\varPsi}

\renewcommand{\Omega}{\varOmega}

\renewcommand{\inlinecomment}[3]{%
{\color{#1}[#2: #3]}} 
\renewcommand\commenter[2]%
{%
\expandafter\newcommand\csname i#1\endcsname[1]{\inlinecomment{#2}{#1}{##1}} 
\expandafter\newcommand\csname #1\endcsname[1]{\mycomment{#1}{##1}{#2}} 
\expandafter\newcommand\csname ri#1\endcsname[1]{} 
\expandafter\newcommand\csname r#1\endcsname[1]{} 
}

\providecommand\ie{\emph{i.e.}}

\providecommand\eg{\emph{e.g.}}
\providecommand\st{\textit{s.t. }}

\providecommand\cf{\textit{cf. }}

\usepackage{xspace}
\usepackage{bm}
\newcommand{\algname}[1]{{\sf\footnotesize#1}\xspace}

\usepackage[capitalize,noabbrev]{cleveref} 
\crefname{relctr}{relation}{Relation} 

\usepackage{pifont}
\usepackage{etoolbox}
\newcounter{relctr} 
\everydisplay\expandafter{\the\everydisplay\setcounter{relctr}{0}} 
\makeatletter
\def\blackcirclenumbers#1{\expandafter\@blackcirclenumbers\csname c@#1\endcsname}
\def\@blackcirclenumbers#1{%
\ifcase#1\or \ding{182}\or \ding{183}\or \ding{184}\or \ding{185}\or \ding{186}\or \ding{187}\or \ding{188}\or \ding{189}\or \ding{190}\or \ding{191}\or \else\@ctrerr\fi}
\makeatother
\AtBeginDocument{} 
\usepackage{tcolorbox}
\usepackage{pifont}
\definecolor{mydarkgreen}{RGB}{39,130,67}
\definecolor{mydarkred}{RGB}{192,25,25}
\newcommand{\green}{\color{mydarkgreen}}
\newcommand{\red}{\color{mydarkred}}
\newcommand{\cmark}{\green\ding{51}}%
\newcommand{\xmark}{\red\ding{55}}%

\let\oldzeta\zeta
\renewcommand{\zeta}{{\green \oldzeta}}
\providecommand{\barzeta}{{\red \bar \oldzeta}}
\let\olddelta\delta
\renewcommand{\delta}{{\green \olddelta}}
\providecommand{\bardelta}{{\red \bar \olddelta}}
\let\oldrho\rho
\renewcommand{\rho}{{\green \oldrho}}
\let\oldcM\cM
\renewcommand{\cM}{{\green \oldcM}}
\providecommand{\barL}{{\red \bar L}}

\providecommand{\FS}{{\green fully stochastic}}
\providecommand{\SMP}{{\red stochastic multi-point}}

\usepackage[round]{natbib}

\begin{document}

%
\runningtitle{Revisiting LocalSGD and SCAFFOLD}

%
\runningauthor{Luo, Stich, Horv\'{a}th, Tak\'{a}\v{c}}

\twocolumn[

\aistatstitle{Revisiting LocalSGD and SCAFFOLD: \\Improved Rates and Missing Analysis}

\aistatsauthor{ Ruichen Luo \And Sebastian U. Stich \And  Samuel Horv\'{a}th \And Martin Tak\'{a}\v{c} }

\aistatsaddress{ IST Austria \\ rluo@ist.ac.at \And CISPA Helmholtz Center \\ stich@cispa.de \And MBZUAI \\ samuel.horvath@mbzuai.ac.ae \And MBZUAI \\ Takac.MT@gmail.com  }
]

\begin{abstract}

\algname{LocalSGD} and \algname{SCAFFOLD} are widely used methods in distributed stochastic optimization, with numerous applications in machine learning, large-scale data processing, and federated learning. However, rigorously establishing their theoretical advantages over simpler methods, such as minibatch SGD (\algname{MbSGD}), has proven challenging, as existing analyses often rely on strong assumptions, unrealistic premises, or overly restrictive scenarios.

In this work, we revisit the convergence properties of \algname{LocalSGD} and \algname{SCAFFOLD} under a variety of existing or weaker conditions, including gradient similarity, Hessian similarity, weak convexity, and Lipschitz continuity of the Hessian. Our analysis shows that (i) \algname{LocalSGD} achieves faster convergence compared to \algname{MbSGD} for weakly convex functions without requiring stronger gradient similarity assumptions; (ii) \algname{LocalSGD} benefits significantly from higher-order similarity and smoothness; and (iii) \algname{SCAFFOLD} demonstrates faster convergence than \algname{MbSGD} for a broader class of non-quadratic functions. These theoretical insights provide a clearer understanding of the conditions under which \algname{LocalSGD} and \algname{SCAFFOLD} outperform \algname{MbSGD}.
\end{abstract}

\section{Introduction}

This paper focuses on solving the following \emph{distributed non-convex optimization} problem: 
\begin{equation}
\label{eq:finite-sum}
\min _{x \in \R^d} \ f(\xx) := \frac 1 n \sum _{i=1} ^{n} f_i (\xx) , 
\end{equation}
where $\xx \in \R^d$ is the optimization variable, $f$ is the global objective function, $n$ is the number of workers, and $f_i$ is the local objective function distributed to the $i$th worker, $i \in [n] := \{ 1, \cdots, n \}$. The local objective functions are \emph{heterogeneous} in general, \ie, $f_i \neq f$. 

\textbf{Distributed non-convex optimization.}
Large scale optimization problems arise in various fields of science and engineering, such as management science~\citep{magnanti2021optimization}, signal processing~\citep{so2020nonconvex}, machine learning~\citep{bottou2018optimization}, among others. 
Given the scale of the problems, it is essential to consider distributed settings, where the optimization tasks are divided among multiple workers who work in parallel to solve the problem. 
In many real-world applications, the objective functions are often non-convex, making the optimization process more challenging compared to convex problems.  

\textbf{Stochastic gradient methods with intermittent communication.} 
For non-convex optimization problems, stochastic gradient descent (\algname{SGD})~\citep{robbins1951sgd} is a basic and strong method widely used, \eg, in training machine learning models~\citep{bottou2010large}. 
When applying \algname{SGD} in distributed settings, the communication between the workers often becomes the bottleneck~\citep{jaggi2014communication,mcmahan2017communication}. 
To adapt \algname{SGD} to distributed setting, many communication-efficient algorithms are designed so that the workers do more local computations and communicate only intermittently.

\textbf{MbSGD vs.\ LocalSGD/SCAFFOLD.}  
Among the communication-efficient methods, we focus on three. Mini-batch SGD (or \algname{MbSGD})~\citep{woodworth2020minibatch,takavc2013mini} is the simplest baseline and is well-understood, while \algname{LocalSGD}~\citep{mcmahan2017communication,stich2018local} and \algname{SCAFFOLD}~\citep{karimireddy2020scaffold} are another two methods gaining great popularity recently. 
At their core, all three methods aim to reduce communication by performing more local computations between communication rounds. 
In \algname{MbSGD}, each worker computes a more accurate stochastic gradient using a larger mini-batch; in \algname{LocalSGD}, each worker takes several local gradient steps before communicating; and in \algname{SCAFFOLD}, those local steps are enhanced with SAGA-like variance reduction techniques~\citep{defazio2014saga}. 
All these methods have been widely adopted in distributed optimization, with their original papers amassing thousands of citations as of 2024. Their importance is reflected in their integration into open-source federated learning frameworks such as \algname{FedML}~\citep{he2020fedml}, \algname{NVIDIA Flare}~\citep{roth2022nvidia}, and \algname{Flower}~\citep{beutel2020flower}.  
While in practice, \algname{LocalSGD} and \algname{SCAFFOLD} often achieve faster convergence than \algname{MbSGD}, there remains a notable gap between theory and practice. 
From a theoretical point of view, \algname{MbSGD} is a formidable baseline, especially in the presence of {heterogeneity}, where it dominates most of the existing analyses.

\textbf{Open questions and missing analysis.} 
Despite the growing interest in distributed non-convex optimization, popular methods like \algname{LocalSGD} and \algname{SCAFFOLD} remain not well understood even under the classic assumptions, such as gradient
similarity (\emph{a.k.a.} bounded gradient dissimilarity), Hessian similarity (\emph{a.k.a.} bounded Hessian dissimilarity), weak convexity, and Lipschitz continuous Hessian.
Several theoretical questions remain open. For instance, is there any speedup analysis of \algname{LocalSGD} for non-convex functions? Can \algname{LocalSGD} converge faster without the assumption of uniform gradient similarity? Can \algname{LocalSGD} benefit from higher-order conditions? Is there any speedup analysis of \algname{SCAFFOLD} for general non-convex functions, potentially without the stringent assumption of uniform Hessian similarity? And to what extent can we relax the assumptions for quadratic functions? 
Given the growing body of work in this field and the significant impact of these two algorithms, it is crucial to properly understand their convergence behaviors, at least under the classic assumptions. Moreover, since these two algorithms have become standard baselines for communication-efficient distributed optimization, ensuring that their convergence rates are cited correctly is essential for fair algorithmic comparisons and future advancements in the field. 

\textbf{Contributions.}
To address the above issues, we revisit the convergence rates of \algname{LocalSGD} and \algname{SCAFFOLD} under gradient similarity, Hessian similarity, weak convexity, and Lipschitz Hessian. 
Our key contributions are as follows: 

\begin{compactitem}
    \item We show that \algname{LocalSGD} can converge faster than \algname{MbSGD} for weakly convex functions, marking the first time such a speedup has been extended to non-convex analysis. 
    \item We show that standard gradient similarity is sufficient for \algname{LocalSGD} to achieve a speedup, making the stronger condition of uniform gradient similarity (as used in~\cite{woodworth2020minibatch,patel2024limits}) redundant. 
    \item We show that \algname{LocalSGD}, a basic gradient method without variance reduction or quasi-Newton tricks, can also benefit from higher-order conditions. 
    \item We show that \algname{SCAFFOLD} can converge faster than \algname{MbSGD} under standard Hessian similarity and weak convexity. This is the first analysis demonstrating such a speedup for general non-quadratic functions without relying on uniform Hessian similarity. 
    \item We introduce a weaker assumption: the existence of a Lipschitz continuous function within the convex hull of the distributed functions. This assumption underpins our proof that \algname{LocalSGD} benefits from higher-order conditions and also allows us to obtain a similar speedup result for \algname{SCAFFOLD} as seen in~\cite{karimireddy2020scaffold}. 
\end{compactitem}

We summarize our results and compare them with the existing analyses of different methods in~\cref{tab:rates}. 

\begin{table*}[ht]
\caption{\small{Summary of the asymptotic non-convex rates of the relevant gradient-based distributed optimization algorithms. Parameters include $\zeta$ (or $\barzeta$) for gradient similarity, $\delta$ (or $\bardelta$) for Hessian similarity, $\rho$ for weak convexity, and $\cM$ for the Lipschitz continuity of the Hessian. ``Speedup'' refers to whether the algorithm can converge faster than the \algname{MbSGD} baseline. ``General'' indicates there is no additional restriction on the class of functions. ``Deterministic'' specifies if the algorithm's randomness only comes from the stochastic oracle. ``Oracle'' refers to the type of gradient oracle used in the algorithm. Detailed definitions are provided in~\cref{sec:problem definition}. 
}} 
\label{tab:rates}
\resizebox{\textwidth}{!}
{\begin{minipage}{1.80\textwidth}
\renewcommand{\footnotesize}{\large}
\centering
\begin{tabular}{@{}cccccccc@{}}
\toprule
{\textbf{Algorithm}} & \textbf{Analysis} & \textbf{Suboptimality} & \textbf{Speedup} & \textbf{General} & \textbf{Deterministic} & \textbf{Oracle} 
\\
\midrule 
\vspace{1mm}
\algname{MbSGD}
&{{\cite{dekel2012optimal}}} 
& $\displaystyle \frac {{ L}\Delta} {R} + \sqrt {\frac {L\Delta \sigma^2} {n\tau R}} $ 
\footnote{The suboptimality of $\frac 1 T \sum_{t=0} ^{T-1} \norm{\nabla f(\bar \xx_t)}^2$. This is standard in the analysis of non-convex functions~\citep{nesterov2003introductory}. \label{ft:average non-convex}}
& baseline 
& \cmark 
& \cmark
& \FS

\vspace{1mm}
\\
\midrule
\vspace{1mm}
\rule{0pt}{4ex} 
\multirow{2}{*}[-.75em]{
    \algname{LocalSGD} 
}
& {{\cite{koloskova2020unified}}} 
& $\displaystyle \frac { { L} \Delta} {R} + \sqrt {\frac {L\Delta \sigma^2} {n\tau R}} + \left( \frac { { L}\Delta { \zeta}} {R} \right)^{\frac23} + \frac {(L\Delta\sigma)^{\frac23}} {\tau^{\frac13}R^{\frac23}}$ 
\footref{ft:average non-convex}
& \xmark 
& \cmark
& \cmark
& \FS

\\
\rule{0pt}{4ex} 

& {{\cite{woodworth2020minibatch}}}
& $\displaystyle \frac {LD^2} {{ \tau} R} + \frac {\sigma D} {\sqrt{n \tau R}} + \left( \frac {L { \barzeta}^2 D^4} {R^2} \right) ^{\frac 1 3} + \left( \frac {L \sigma^2 D^4} {\tau R^2} \right) ^{\frac 1 3}$
\footnote{The suboptimality of $\frac 1 T \sum _{t=0} ^{T-1} \Eb {f(\bar \xx_t)} - f^*$. This is standard in the analysis of convex functions~\citep{nesterov2003introductory}.  \label{ft:average convex}}
& \cmark
& {\red convex}
& \cmark
& \FS

\vspace{1mm}
\\
\midrule
\vspace{1mm}
\rule{0pt}{4ex} 
\multirow{2}{*}[-.75em]{
    \algname{SCAFFOLD} 
}
& {{{\cite{karimireddy2020scaffold}}}}
& $\displaystyle \frac { { L} \Delta} {R} + \sqrt {\frac {L\Delta \sigma^2} {n\tau R}}$ 
\footnote{The suboptimality of $\frac 1 R \sum _{r=0} ^{R-1} \E \norm { \nabla f(\bar \xx_{2r \tau}) }^2$. In their analysis, they combine all the iterates between two communication rounds into one big step and analyze the descent. Moreover, for theoretical purposes, they use a different global stepsize $\eta_g$ to aggregate the workers' updates. \label{ft:average at comm round}}
& \xmark
& \cmark
& \cmark
& \FS

\\
\rule{0pt}{4ex}

& {{{\cite{karimireddy2020scaffold}}}}
& $\displaystyle \left( { \frac L \tau + \bardelta + \rho} \right) \frac {\Delta} {R} + \sqrt {\frac {L\Delta \sigma^2} {n\tau R}}$ 
\footref{ft:average non-convex}
& \cmark
& {\red quadratic}
& \cmark
& {\red exact}
\footnote{They use exact gradients in the variance reduction term.}

\vspace{1mm}
\\
\midrule
\vspace{1mm}
\rule{0pt}{4ex} 
\algname{ScaffNew}
& {{\cite{mishchenko2022proxskip}}}
&
$\displaystyle (1-\mu/L)^T D^2 + (1-\theta)^T D^2 + \eta^2 \frac {\sigma^2} {\theta} $
\footnote{The suboptimality of $\E \norm {\bar \xx_T - \xx^*}^2$. This is standard in the analysis of strongly-convex functions~\citep{nesterov2003introductory}. In this rate, $\theta := \min \{ \mu \eta, p^2 \}$, where $p$ is the communication probability. \label{ft:last iterate strong convexity}}
&
\cmark
&
{\red $\mu$-strongly convex}
&
\xmark
&
\FS

\vspace{1mm}
\\
\midrule
\vspace{1mm}
\rule{0pt}{4ex} 
\algname{CE-LSGD}
& {{\cite{patel2022towards}}}
&
$\displaystyle \left( { \frac {\barL} {\sqrt{\tau}} + \bardelta} \right) \frac {\Delta} {R} + \frac {{\sigma}^2} {n \tau R} + \left( \frac {\barL \Delta {\sigma}} {n \tau R} \right)^{\frac 2 3}$
\footref{ft:average non-convex}
&
\cmark
&
\cmark
&
\xmark
&
\SMP
\footnote{In \algname{CE-LSGD}, they use a stochastic multi-point gradient oracle, where each worker can query $G_i(\cdot, \xi_t^i)$ again at a different point in the next iteration. Moreover, they need to assume a stronger $\barL$-mean smoothness in their analysis, \ie, $\forall \xx, \yy \in \R^d$, $\E _{\xi_t^i} \norm {G_i(\xx, \xi_t^i) - G_i(\yy, \xi_t^i)} \le \barL \norm {\xx - \yy}$.}

\vspace{1mm}
\\
\midrule
\vspace{1mm}
\rule{0pt}{4ex} 
\multirow{3}{*}[-1.5em]{\algname{LocalSGD}}
&
    \textbf{{\cref{thm:LocalSGD faster}}} 

& $\displaystyle \left( { \frac {L} {\tau} + \rho} \right) \frac {\Delta} {R} + \sqrt {\frac {L\Delta \sigma^2} {n\tau R}} + \left( \frac {L\Delta { \zeta}} {R} \right)^{\frac 2 3} + \frac {(L\Delta\sigma)^{\frac 2 3}} {\tau^{\frac 1 3}R^{\frac 2 3}}$ 
\footref{ft:average non-convex}
& \cmark
& \cmark
& \cmark
& \FS

\\
\rule{0pt}{4ex} 

& \textbf{{\cref{thm:LocalSGD faster convex}}} 
& $\displaystyle \frac {LD^2} {{ \tau} R} + \frac {\sigma D} {\sqrt{n \tau R}} + \left( \frac {L { \zeta}^2 D^4} {R^2} \right) ^{\frac 1 3} + \left( \frac {L \sigma^2 D^4} {\tau R^2} \right) ^{\frac 1 3}$ 
\footref{ft:average convex}
& \cmark 
& \red{convex}
& \cmark
& \FS

\\
\rule{0pt}{4ex} 

& \textbf{{\cref{thm:LocalSGD faster} \& \ref{thm:LocalSGD improved conditioning}}} 
\footnote{This rate can be obtained by combining the proofs of \cref{thm:LocalSGD faster,thm:LocalSGD improved conditioning}.} 
& $\displaystyle \left( { \frac {L} {\tau} + \rho} \right) \frac {\Delta} {R} + \sqrt {\frac {L\Delta \sigma^2} {n\tau R}} + \left( \frac {\bardelta \Delta { \zeta}} {R} \right)^{\frac 2 3} + \frac {(L\Delta\sigma)^{\frac 2 3}} {\tau^{\frac 1 3}R^{\frac 2 3}} + \left( \frac {\cM^2\Delta^4\zeta^4} {R^4} \right)^{\frac15}$
\footref{ft:average non-convex}
& \cmark
& \cmark
& \cmark
& \FS

\vspace{1mm}
\\
\midrule
\vspace{1mm}
\rule{0pt}{4ex} 
\multirow{2}{*}[-.75em]{\algname{SCAFFOLD}}

& \textbf{\cref{thm:SCAFFOLD speedup}} 
& $\displaystyle \left( { \frac L \tau + \sqrt{L\delta} + \rho} \right) \frac {\Delta} {R} + \sqrt { \frac {L \Delta \sigma^2} {n \tau R} } + \frac {(L\Delta\sigma)^{\frac 2 3}} {\tau^{\frac 1 3}R^{\frac 2 3}} $ 
\footref{ft:average non-convex}
& \cmark
& \cmark
& \cmark
& \FS

\\
& 
    \textbf{\cref{thm:Lipschitz Hessian SCAFFOLD}} 

& $\displaystyle \left( { \frac L \tau + \sqrt {\bardelta \delta} + \rho} \right) \frac {\Delta} {R} + \sqrt { \frac {L \Delta \sigma^2} {n \tau R} } + \frac {(\bardelta\Delta\sigma)^{\frac 2 3}} {\tau^{\frac 1 3}R^{\frac 2 3}} $ 
\footref{ft:average non-convex}
& \cmark
& {\red \cref{assumption:Lipschitz Hessian} ($\cM=0$)}
& \cmark
& \FS

\vspace{1mm}
\\
 
\bottomrule
\end{tabular}
\end{minipage}}
\end{table*}

\textbf{Structure of the paper.} 
We discuss the related work in~\cref{sec:work}. Then, we formally define the problem, algorithms, and assumptions in ~\cref{sec:problem definition}. We review the current rates of \algname{MbSGD}, \algname{LocalSGD}, and \algname{SCAFFOLD} in~\cref{sec:existing analysis}. Our new analysis is presented in~\cref{sec:analysis}, with key implications and technical innovations highlighted in the main text, while detailed proofs are deferred to~\cref{appendix_sec:proof_details}. In~\cref{sec:expr}, we present suitable synthetic experiments that validate our theoretical findings. We conclude with a discussion of limitations and future work.

\section{Related Work}
\label{sec:work}

\textbf{Classic assumptions and existing analysis.} 
Early analysis focuses on the more restrictive, homogeneous setting~\citep{zhou2018on,yu2019parallel,stich2018local,woodworth2020local}. This paper, instead, considers the more general, heterogeneous settings. Recent works on the heterogeneous settings often consider the following classic assumptions: gradient similarity, Hessian similarity, weak convexity, and Lipschitz continuous Hessian. 
Below, we briefly summarize the existing analyses based on these assumptions: 
\begin{compactitem}
    \item \cite{koloskova2020unified,karimireddy2020scaffold} show that \algname{LocalSGD} can only match \algname{MbSGD} under standard gradient similarity. Moreover, the speedup of \algname{LocalSGD} is only proven for convex functions, and this relies on a stronger condition of uniform gradient similarity~\citep{woodworth2020minibatch}. 
    \item For quadratic functions, \cite{karimireddy2020scaffold} show the speedup of \algname{SCAFFOLD} under uniform Hessian similarity and weak convexity. However, for more general non-quadratic functions, there is no theoretical proof of a speedup for \algname{SCAFFOLD}. 
\end{compactitem}

\textbf{Related work on new algorithms.}
Distributed non-convex optimization continues to be an active research area. 
Perhaps motivated by the limitations of the current theories of \algname{LocalSGD} and \algname{SCAFFOLD}, many recent works focus on designing new algorithms. Examples include \algname{ProxSkip} (or \algname{ScaffNew})~\citep{mishchenko2022proxskip}, \algname{BVR-L-SGD}~\citep{murata2021bias}, \algname{CE-LSGD}~\citep{patel2022towards}, \algname{MimeMVR}~\citep{karimireddy2021breaking}, \algname{DANE+} and \algname{FedRed}~\citep{jiang2024federated}. Each of these works considers somehow different settings. \cite{mishchenko2022proxskip} proposes \algname{ScaffNew}, a variant of the \algname{SCAFFOLD} algorithm with randomized communication rounds, and shows acceleration from local steps, but their result holds only in expectation for strongly convex functions. \cite{murata2021bias,patel2022towards} study algorithms with stochastic multi-point oracle and randomized worker sampling, while the algorithms considered in this paper are deterministic and only use the fully stochastic gradient oracle, which is more general. \cite{karimireddy2021breaking} assumes each local function $f_i$ is a finite sum of components with similar Hessians across workers, making the setting more restrictive. \cite{jiang2024federated} studies variance-reduced proximal methods, which are more difficult to implement, and it is not clear how to handle stochastic noise. 
Thus far, we are not aware of any algorithm that uses a fully stochastic oracle and converges provably better than \algname{MbSGD}/\algname{LocalSGD}/\algname{SCAFFOLD} for general non-convex functions. 

\textbf{Related work on new assumptions.}
A parallel line of research focuses on analyzing \algname{LocalSGD} under new assumptions. \cite{wang2022unreasonable} introduces a novel assumption that the (so-called) average pseudo-gradient is bounded at optimum. However, this analysis is problematic: if the assumed bound is near zero, \algname{LocalSGD} converges in constant communication rounds; otherwise, it diverges. Moreover, \cite{patel2023on} proves the lower bound, showing that this assumption fails for non-convex functions. \cite{zindari2023on} considers quadratic functions under the assumption of identical Hessians or identical optima, but these are overly restrictive and do not apply to general, non-convex functions. Ultimately, none of these works provide a reasonable explanation for the fast convergence observed in practice for \algname{LocalSGD}. 

\textbf{Related work on generalization.}
In the context of distributed learning, our optimization problem~\eqref{eq:finite-sum} corresponds to empirical risk minimization. Some more recent works, under different set of assumptions, focus on bounding the generalization gap~\citep{sefidgaran2024lessons,sun2024understanding,gu2023why}. These works are orthogonal or somewhat complementary to this paper, as the true risk in distributed learning is bounded by the sum of the empirical risk and the generalization gap. We include a more detailed discussion in~\cref{appendix_sec:generalization}.  

\section{Preliminary} 
\label{sec:problem definition} 

We start with the definition of the optimization problem. Then, we describe the algorithms of interest and, finally, the technical assumptions. 

\subsection{Distributed Stochastic Non-Convex Optimization with Intermittent Communication}
\label{sec:setup}

\textbf{Distributed non-convex optimization.}
We recall the distributed non-convex optimization problem: 
\[
\min _{x \in \R^d} \ f(\xx) := \frac 1 n \sum _{i=1} ^{n} f_i (\xx) , 
\]
and we assume throughout the paper that $f$ is bounded below, and $f_i \colon \R^d \rightarrow \R$ is possibly non-convex and has $L$-Lipschitz continuous gradient ($L > 0$), for all $i \in [n]$. 

\textbf{Stochastic oracle.}
We focus on solving problem~\eqref{eq:finite-sum} by iterative algorithms via subsequent queries to a \emph{fully stochastic oracle} $\mathcal{SO}$. Let $T$ be the total number of successive queries to $\mathcal{SO}$. 
At iteration $t$, where $0 \le t \le T-1$, the workers' inputs being $\left( \xx_t^1, \cdots, \xx_t^n \right) \in \R^{d\times n}$, the $\mathcal {SO}$ outputs vectors 
\[
  \left( \gg_t^1, \cdots, \gg_t^n \right) := \left( G_1(\xx_t^1, \xi_t^1), \cdots, G_n(\xx_t^n, \xi_t^n) \right) \in \R^{d\times n},
\]
where $\left\{\xi_t^i : {0 \le t \le T-1} \right\}$ are i.i.d. random variables for each $i \in [n]$. 
We assume the following conditions on the Borel functions $G_i(\xx, \xi_t^i)$: 
\begin{equation}
\begin{gathered}
\label{eq:borel}
    \EEb {\xi_t^i} {G_i (\xx, \xi_t^i)} = \nabla f_i(\xx), \\ 
    \EE {\xi_t^i} {\norm {G_i (\xx, \xi_t^i) - \nabla f_i(\xx)}^2} \le \sigma^2 .
\end{gathered}
\end{equation}

\textbf{Intermittent communication.}
In distributed optimization, each worker $i \in [n]$ has access only to its local function $f_i$, which is not visible to the other workers. 
The optimization process proceeds through \emph{intermittent communication}, where workers communicate periodically after every \emph{interval} of $\tau \in \mathbb Z _{\ge 2}$ iterations. We assume w.l.o.g.\ that $T$ is a multiple of $\tau$. That is, after the $r$-th communication round ($r \in [0, T/\tau-1]$), each worker $i \in [n]$ can access the gradients $\{ \gg_{r\tau}^i, \cdots, \gg_{r\tau+\tau-1}^i \}$, and only then communicates with the other workers at the subsequent $(r+1)$-th communication round. 

To simplify the notations, let $\bar \xx_t = \frac 1 n \sum _{i=1} ^n \xx_t^i$, $f^* = \inf _{\xx \in \R^d} f(\xx)$, and $\Delta = f(\bar \xx_0) - f^*$. 

\subsection{Algorithms} 

In the context of distributed non-convex optimization with intermittent communication, we focus on the analysis of three algorithms, \algname{MbSGD}, \algname{LocalSGD}, and \algname{SCAFFOLD}. Below, we give their formal descriptions. 
In these algorithms, we assume $\xx_0^1 = \cdots = \xx_0^n$.

\textbf{{MbSGD} and {LocalSGD}.}
Both \algname{MbSGD} and \algname{LocalSGD} are two well-recognized, basic methods in distributed optimization. Let $R$ be the total number of communication rounds, \ie, $R = T/ \tau$.  
\begin{compactitem}
    \item \textbf{{MbSGD.}} \ Between two communication rounds, each worker computes its gradient estimate on a $\tau$-times larger mini-batch and communicates at the next communication round. Formally, for each $i \in [n]$ and $t \in [0,T-1]$, 
    \begin{equation}
      \xx _{t+1} ^i = \begin{cases}
        \bar \xx_{t-\tau+1} - \frac {\eta} {n} \sum _{j=1} ^n \sum _{k=0}^{\tau-1} \gg_{t-k}^j, \qquad \qquad & \\
        \hfill \text{ if } t+1 \text{ is a multiple of } \tau, & \\
        \xx_t^i, \hfill \text{ otherwise}. & 
      \end{cases}
    \end{equation}
    
    \item \textbf{{LocalSGD.}} \ Between two communication rounds, each worker takes $\tau$ local steps using its gradient estimates. Then, at the next communication round, all workers average their solutions globally. Formally, for each $i \in [n]$ and $t \in [0,T-1]$, 
    \begin{equation}
    \label{eq:local sgd}
      \xx _{t+1} ^i = \begin{cases}
        \bar \xx_{t-\tau+1} - \frac {\eta} {n} \sum _{j=1} ^n \sum _{k=0}^{\tau-1} \gg_{t-k}^j, \qquad \qquad & \\ 
        \hfill \text{ if } t+1 \text{ is a multiple of } \tau, & \\ 
        \xx_t^i - \eta \gg_t^i, \hfill \text{ otherwise}. &
      \end{cases}
    \end{equation}
\end{compactitem} 

\textbf{SCAFFOLD.}
\algname{SCAFFOLD}, introduced by~\citet{karimireddy2020scaffold}, incorporates variance reduction into the local steps to address the issue of client drift. Suppose there are $2R$ communication rounds, \ie, $2R = T/\tau$. 
The algorithm considered in this paper is a simplified variant of the original \algname{SCAFFOLD} algorithm.
As described in~\cref{alg:SCAFFOLD}, across two consecutive communication rounds, in the first round, each worker computes a stochastic gradient on a $\tau$-times larger mini-batch; and in the second round, they take variance-reduced gradient steps using the gradients computed in the first round. 
The original algorithm in~\cite{karimireddy2020scaffold} is a bit more complex, as it allows partial worker participation, the use of the last update for variance reduction, and different local/global stepsizes. We drop these features because none of them are supported in their speedup analysis (\cf Section~6 of \cite{karimireddy2020scaffold}). Also, the use of a different global stepsize is mainly for theoretical purposes and is not a common practice in real-world applications (\cf Section~7.1 of \cite{karimireddy2020scaffold}). Hence, our simplified version retains the core variance reduction mechanism while streamlining the analysis.

\begin{algorithm}[ht!]
\caption{SCAFFOLD}
\begin{small}
\begin{algorithmic}[1]
    \For {$r=0,1,\cdots, R-1$}
        \For {$i \in [n]$} \textbf{in parallel}
        \For {$k=0,1,\cdots, \tau-1$}
            \State $\xx_{2r \tau + k + 1}^i = \xx_{2r \tau + k}^i$
        \EndFor
        \State $\hat \gg _{(r \tau)} ^i = \frac 1 \tau \sum _{k=0} ^{\tau - 1} \gg _{2r \tau + k} ^i $ \label{line:variance reduction in scaffold}
        \EndFor

        \State Compute and broadcast: \(\hat \gg _{(r \tau)} = \frac 1 n \sum _{i=1} ^n \hat \gg _{(r \tau)} ^i\)

        \For {$i \in [n]$} \textbf{in parallel}
        \For {$k=\tau,\tau+1,\cdots, 2\tau-2$}
            \State \vspace{-1em}
            \begin{equation} 
            \begin{aligned}
            \label{eq:scaffold_update}
            \qquad \quad \xx_{2r \tau + k + 1}^i &= \xx_{2r \tau + k}^i  - \eta \left( \gg _{2r \tau + k} ^i - \hat \gg _{(r \tau)} ^i + \hat \gg_{(r \tau)} \right) 
            \end{aligned}
            \end{equation}
            \vspace{-1em}
        \EndFor
        \EndFor

        \State Compute: 
        \(
        \bar \xx_{2(r+1)\tau} = \bar \xx _{2r \tau} - \frac \eta n \sum _{j=1} ^n \sum _{l=\tau} ^{2\tau-1} \gg_{2r \tau + l}^j
        \)
        \label{line:aggregation in scaffold}
        \State Broadcast: 
        \(
            \xx_{2(r+1)\tau}^i = \bar \xx_{2(r+1)\tau} , \text{ for each } i \in [n]
        \)
    \EndFor
\end{algorithmic}
\end{small}
\label{alg:SCAFFOLD}
\end{algorithm} 

\subsection{Assumptions}

We present the technical assumptions in this section and discuss their relevance and implications. 

Standard gradient similarity is required to prove the convergence of \algname{LocalSGD}. Similar assumptions are used in many previous works~\citep{koloskova2020unified,karimireddy2020scaffold}. 

\begin{assumption}[Standard gradient similarity]
\label{assumption:GS}
    For some ${ \zeta} \ge 0$, we have  
    \begin{equation}
    \label{eq:GS bound}
      \sup _{\xx \in \R^d} \frac 1 n \sum _{i=1}^n \norm {\nabla f_i(\xx) - \nabla f(\xx)}^2 \le { \zeta}^2 . 
    \end{equation}
\end{assumption}

Some recent works, such as \cite{woodworth2020minibatch,patel2024limits}, rely on the stronger assumption of uniform gradient similarity. 
\begin{strongerassumption}[Uniform gradient similarity]
\label{assumption:UGS}
For some $\barzeta \ge 0$, we have
\begin{equation}
    \label{eq:uniform GS}
    \sup _{\xx \in \R^d} \sup _{i \in [n]} \norm {\nabla f_i(\xx) - \nabla f(\xx)}^2 \le { \barzeta^2} 
\end{equation}
\end{strongerassumption}

\begin{remark}
    Based on the definitions, $\zeta$ is always less than or equal to $\barzeta$. In cases where the heterogeneity is very imbalanced between the workers, $\zeta$ could be significantly smaller--up to a factor of $n^{1/2}$. Our speedup analysis of \algname{LocalSGD} avoids the need for this stronger assumption, which makes it applicable to scenarios with more heterogeneous data.  
\end{remark}

To study the speedup of taking local steps in \algname{LocalSGD} and \algname{SCAFFOLD}, following~\cite{karimireddy2020scaffold,karimireddy2021breaking}, we introduce the assumptions of standard/uniform Hessian similarity and weak convexity: 

\begin{assumption}[Standard Hessian similarity]
\label{assumption:mean HS}
    For some $\delta \in [0,L]$, we have 
    \begin{equation}
    \label{eq:mean HS bound}
    \begin{aligned}
    &\frac 1 n \sum_{i=1} ^n \norm { \nabla f_i(\xx) - \nabla f(\xx) - \nabla f_i(\yy) + \nabla f(\yy) }^2 \\
    &\qquad \le {\delta^2} \norm {\xx - \yy}^2 ,
    \end{aligned}
    \end{equation}
    for all $\xx, \yy \in \R^d$. 
\end{assumption}

\begin{strongerassumption}[Uniform Hessian similarity]
\label{assumption:HS}
    For some $\bardelta \in [0,2L]$, we have   
    \begin{equation}
    \label{eq:HS bound}
    \norm { \nabla f_i(\xx) - \nabla f(\xx) - \nabla f_i(\yy) + \nabla f(\yy) } \le {\bardelta} \norm {\xx - \yy} ,
    \end{equation}
    for all $\xx, \yy \in \R^d$ and for all $i \in [n]$. 
\end{strongerassumption}

\begin{assumption}[Weak convexity]
\label{assumption:weak convexity}
    For some $\rho \in [0,L]$, we have 
    \[
    f_i(\xx) + \frac {\rho} 2 {\xx \T \xx} \text{ is convex} ,
    \]
    for all $i \in [n]$. 
\end{assumption}

\begin{remark}
    We remark on the above assumptions. 
    \begin{compactitem}
    \item
        Hessian similarity has been widely used in the literature of distributed optimization~\citep{shamir2014communication,karimireddy2020scaffold,karimireddy2021breaking}. Despite the name, we do not necessitate that $f_i$ is twice continuously differentiable. Indeed, all functions with $L$-Lipschitz continuous gradient satisfy $L$-standard Hessian similarity and $2L$-uniform Hessian similarity. 
    \item
        Weak convexity is often used in distributed optimization to show the effectiveness of the local gradient methods~\citep{karimireddy2020scaffold,karimireddy2021breaking}. As a matter of fact, all the functions with $L$-Lipschitz continuous gradients are $L$-weakly convex, and the convex functions are $0$-weakly convex. 
    \end{compactitem} 
    In this paper, when adopting these assumptions, we always assume $\delta \lll L$, $\bardelta \lll 2L$, and $\rho \lll L$. 
\end{remark}

Finally, we introduce a new variant of Lipschitz continuous Hessian.  

\begin{assumption}[Lipschitz continuous Hessian]
\label{assumption:Lipschitz Hessian}
    For some $\cM \ge 0$, there exists (at least) one function $\hat f$ such that: 
    \(
      \hat f \in \mathbf{conv} \{f_1, \cdots, f_n\} 
    \), \footnote{ The convex hull \(\mathbf{conv} \{f_1, \cdots, f_n\}\) denotes \(\left\{ \sum _{i \in [n]} a_i f_i \mid \sum _{i \in [n]} a_i = 1 \text{ and } a_i \ge 0, \text{ for all } i \in [n] \right\}\).}
    and   
    \begin{equation}
    \label{eq:Lipschitz Hessian}
    \norm {\nabla^2 \hat f(\xx) - \nabla^2 \hat f(\yy)} \le \cM \norm {\xx - \yy} , 
    \end{equation}
    for all $\xx, \yy \in \R^d$.
\end{assumption}

\begin{remark}
    Quadratic or Lipschitz continuous Hessian has been widely assumed in the literature of distributed optimization~\citep{shamir2014communication,hendrikx2020statistically,karimireddy2020scaffold}. Nonetheless, our \cref{assumption:Lipschitz Hessian} might represent the weakest assumption among similar ones made in the literature. It is significantly weaker than requiring all $f_i$ functions to be quadratic, \ie, $0$-Lipschitz continuous Hessian. Moreover, it is also weaker than the assumption that ``there exists $i \in [n]$ such that
    $\nabla^2 f_i(x)$ is $\cM$-Lipschitz continuous.'' \footnote{For instance, for $f_1(x) := \frac 1 2 (x+2)^2 + 2 \log (1+x^2)$ and $f_2(x) := \frac 1 2 (x-1)^2 - \log (1+x^2)$, while $f_1$ has $6$-Lipschitz continuous Hessian and $f_2$ has $3$-Lipschitz continuous Hessian, we have $\hat f (x) = \frac 1 3 f_1 (x) + \frac 2 3 f_2 (x)$ with $0$-Lipschitz continuous Hessian.}
\end{remark}

\section{Existing Convergence Analyses}
\label{sec:existing analysis}

In this section, we formally state the preliminary results of \algname{MbSGD}, \algname{LocalSGD} and \algname{SCAFFOLD} in literature, and discuss the open questions and missing analysis. 

We first state the well known rate of \algname{MbSGD}. This rate remains tight under all of \cref{assumption:GS,assumption:UGS,assumption:mean HS,assumption:HS,assumption:weak convexity,assumption:Lipschitz Hessian}. 
\begin{lemma}[{\cite{dekel2012optimal}}]
\label{lemma:MbSGD}
    There exists $\eta > 0$ such that \algname{MbSGD} ensures the following upper bound on $\frac 1 T \sum _{t=0}^{T-1} \E \norm {\nabla f(\bar \xx_t)}^2$: 
    \begin{equation}
    \label{eq:MbSGD rate}
    \cO \left( \frac {{ L}\Delta} {R} + \sqrt {\frac {L\Delta \sigma^2} {n\tau R}} \right) .
    \end{equation}
\end{lemma}

\begin{remark}
    As a simplification of \eqref{eq:MbSGD rate}, for the case $L=\Delta=\sigma=1$ and $\tau=\infty$, we get the asymptotic convergence rate of $\cO \left( \frac 1 R \right)$ for \algname{MbSGD}. 
\end{remark}

Now, we state the best-known rate of \algname{LocalSGD} for non-convex functions. 
\begin{lemma}[{\cite{koloskova2020unified}}]
\label{lemma:classic LocalSGD}
    Under \cref{assumption:GS}, there exists $\eta>0$ such that \algname{LocalSGD} ensures the following upper bound on $\frac 1 T \sum _{t=0}^{T-1} \E \norm {\nabla f(\bar \xx_t)}^2$: 
    \begin{equation}
    \label{eq:classic LocalSGD rate}
    \cO \left( \frac { { L} \Delta} {R} + \sqrt {\frac {L\Delta \sigma^2} {n\tau R}} + \left( \frac { { L}\Delta { \zeta}} {R} \right)^{\frac23} + \frac {(L\Delta\sigma)^{\frac23}} {\tau^{\frac13}R^{\frac23}} \right) . 
    \end{equation}
\end{lemma} 

\begin{remark}
\label{remark:classic LocalSGD}
We see that the rate of \algname{LocalSGD} in \eqref{eq:classic LocalSGD rate} is no faster than that of \algname{MbSGD} in \eqref{eq:MbSGD rate}. As a simplification of \eqref{eq:classic LocalSGD rate}, for the case $L=\Delta=\sigma=1$ and $\tau=\infty$, the asymptotic rate of \algname{LocalSGD} is 
\(
\cO \left( \frac { 1} R + \left( \frac {{ \zeta}} {R} \right) ^{2/3} \right) .  
\)
Hence, the conditioning of 
\(
\zeta^2 = \cO \left( 1/R \right)
\) 
is required for \algname{LocalSGD} to match the rate of \algname{MbSGD}. 
\end{remark} 

We also present the only known speedup result of \algname{LocalSGD} for convex functions, as established in \cite{woodworth2020minibatch}. 
\begin{lemma}[\cite{woodworth2020minibatch}]
    \label{lemma:LocalSGD faster convex previous}
    Under \cref{assumption:UGS}, if all the local functions $f_i$ are convex, $\xx^* \in \argmin _{\xx \in \R^d} {f(\xx)}$, and there exists some $D\ge 0$ such that $\norm{ \bar \xx_0 - \xx^*} \le D$, then there exists $\eta > 0$ such that \algname{LocalSGD} ensures the following upper bound on $\frac 1 T \sum _{t=0} ^{T-1} \Eb { f \left( \bar \xx_t \right) } - f^*$: 
    \begin{equation}
    \label{eq:LocalSGD convex faster rate previous}
    \cO \left( \frac {LD^2} {{ \tau} R} + \frac {\sigma D} {\sqrt{n \tau R}} + \left( \frac {L { \barzeta}^2 D^4} {R^2} \right) ^{\frac 1 3} + \left( \frac {L \sigma^2 D^4} {\tau R^2} \right) ^{\frac 1 3} \right) . 
    \end{equation}
\end{lemma}

\begin{remark}
    \Cref{lemma:LocalSGD faster convex previous} shows that \algname{LocalSGD} can converge faster than \algname{MbSGD}. The speedup is reflected in the optimization term, the first term in \eqref{eq:LocalSGD convex faster rate previous}, from $\cO \left( \frac {1} {R} \right)$ to $\cO \left( \frac {1} {\tau R} \right)$. 
    In this speedup analysis, however, the heterogeneity term, the third term in \eqref{eq:LocalSGD convex faster rate previous}, gets worse as compared to \eqref{eq:classic LocalSGD rate}. 
    Suppose we assume only $\zeta$-standard gradient similarity in \cref{lemma:LocalSGD faster convex previous}. Then, the same analysis requires replacing $\barzeta^2$ with $n \zeta^2$ in \eqref{eq:LocalSGD convex faster rate previous}. 
    Therefore, when the heterogeneity is very imbalanced between the workers, this bound can have a dependence of ${ n}^{1/3}$ with respect to the number of workers, which is counterintuitive, as we typically expect faster convergence with the increasing number of workers. 
\end{remark}

\noindent{\bf Open questions for \algname{LocalSGD}.} The following two questions remain open: (i) \emph{whether for non-convex functions there can be a speedup analysis}, and (ii) \emph{whether the stronger condition of uniform gradient similarity can be avoided}. 

Then, we state the convergence rates of \algname{SCAFFOLD}. 

\begin{lemma}[{\cite{karimireddy2020scaffold}}]
\label{lemma:original SCAFFOLD}
    Suppose in \cref{line:aggregation in scaffold} of \cref{alg:SCAFFOLD}, a different global stepsize $\eta_g$ can be used when aggregating the updates of the workers. There exists $\eta_g \ge \eta > 0$ such that \algname{SCAFFOLD} ensures the following upper bound on $\frac 1 R \sum _{r=0}^{R-1} \E \norm {\nabla f(\bar \xx_{2r \tau})}^2$: 
    \begin{equation}
    \label{eq:classic SCAFFOLD rate}
    \cO \left( \frac {{ L}\Delta} {R} + \sqrt {\frac {L\Delta \sigma^2} {n\tau R}} \right) .
    \end{equation} 
\end{lemma}

\begin{lemma}[\cite{karimireddy2020scaffold}]
\label{lemma:original SCAFFOLD quadratic}
    Suppose exact gradients can be used for variance reduction, \ie, $\hat \gg_{(r\tau)} ^i = \nabla f_i(\bar \xx_{2r\tau})$ in \cref{line:variance reduction in scaffold} of \cref{alg:SCAFFOLD}. Under \cref{assumption:HS,assumption:weak convexity}, if all the local functions $f_i$ are quadratic, then there exists $\eta > 0$ such that \algname{SCAFFOLD} ensures the following upper bound on $\frac 2 T \sum _{r=0}^{R-1} \sum _{k=0}^{\tau-1} \E \norm {\nabla f(\bar \xx_{2r \tau + \tau + k})}^2$: 
    \begin{equation}
    \label{eq:quadratic SCAFFOLD rate}
    \cO \left( \left( { \frac L \tau + \bardelta + \rho} \right) \frac {\Delta} {R} + \sqrt {\frac {L\Delta \sigma^2} {n\tau R}} \right) .
    \end{equation} 
\end{lemma}

\begin{remark}
\Cref{lemma:original SCAFFOLD} shows that \algname{SCAFFOLD}, with the rescaling of local and global stepsizes, converges at the same rate as \algname{MbSGD} for general non-convex functions. Further, \cref{lemma:original SCAFFOLD quadratic} shows that, under Hessian similarity and weak convexity, \algname{SCAFFOLD} converges faster than \algname{MbSGD} for quadratic functions. The speedup is reflected in the optimization term, from $\cO \left( \frac 1 {R} \right)$ to $\cO \left( (\frac {1} {\tau} + \bardelta + \rho) \frac {1}{R} \right)$. 
\end{remark}

\noindent{\bf Open question for \algname{SCAFFOLD}.} One key question remains open: 
    \emph{whether any speedup can be proved for general non-quadratic functions and/or without uniform Hessian similarity}. 

Recently, there has been growing interest in designing new algorithms under the various assumptions of Hessian similarity, weak convexity, and/or Lipschitz continuous Hessian. Several works have referred to the rates of \algname{LocalSGD} and/or \algname{SCAFFOLD} as baselines for comparisons, including \cite{karimireddy2020scaffold,karimireddy2021breaking,shen2021stl,khaled2023faster}, to name a few. However, it remains unclear whether these cited rates are tight and correctly interpreted across different settings. Therefore, the following issue is subtle but crucial for fair algorithmic comparisons in the field: 

\noindent{\bf Missing analysis.} The following analysis is missing: \emph{whether the rates of \algname{LocalSGD} and \algname{SCAFFOLD} can be improved under the various classic assumptions}. 

\section{New Convergence Analysis} 
\label{sec:analysis}

To address the aforementioned open problems and missing issues, we revisit the analysis of \algname{LocalSGD} and \algname{SCAFFOLD}. We highlight the technical novelty in our analysis, and defer the proofs to \cref{appendix_sec:proof_details}. 

\subsection{{New Analysis of LocalSGD}} 
\label{sec:localsgd}

We show that the classic rates of \algname{LocalSGD} can be improved in the various settings: (i) Theorem~\ref{thm:LocalSGD faster} shows the speedup under weak convexity and Theorem~\ref{thm:LocalSGD faster convex} shows the improvement for convex functions, both without uniform gradient similarity; (ii) Theorem~\ref{thm:LocalSGD improved conditioning} shows a possible improvement under Hessian similarity and Lipschitz continuous Hessian. 

\begin{theorem}
\label{thm:LocalSGD faster}
    Under \cref{assumption:GS,assumption:weak convexity}, there exists $\eta>0$ such that \algname{LocalSGD} ensures the following upper bound on $\frac 1 T \sum _{t=0}^{T-1} \E \norm {\nabla f(\bar \xx_t)}^2$: 
    \begin{equation}
    \label{eq:LocalSGD faster rate}
    \cO \left( \left( { \frac {L} {\tau} + \rho} \right) \frac {\Delta} {R} + \sqrt {\frac {L\Delta \sigma^2} {n\tau R}} + \left( \frac {L\Delta { \zeta}} {R} \right)^{\frac 2 3} + \frac {(L\Delta\sigma)^{\frac 2 3}} {\tau^{\frac 1 3}R^{\frac 2 3}} \right) . 
    \end{equation}
\end{theorem}

\begin{remark}
As compared to the rate of \algname{MbSGD}~(\cref{lemma:MbSGD}) and the previous rate of \algname{LocalSGD}~(\cref{lemma:classic LocalSGD}), our new rate of \algname{LocalSGD} in \cref{thm:LocalSGD faster} has a faster optimization term of $\cO \left( (\frac L \tau + \rho) \frac \Delta R \right)$. 
As a simplification of \eqref{eq:LocalSGD faster rate}, for the case $L=\Delta=\sigma=1$ and $\tau=\infty$, the asymptotic rate of \algname{LocalSGD} is 
\(
    \cO \left( \frac {{ \rho}} R + \left( \frac {{ \zeta}} {R} \right)^{2/3} \right) .
\)
    We see that the conditioning of 
    \(
    \zeta^2 = \cO \left( 1/R \right) 
    \)
    ensures that \algname{LocalSGD} converges faster than \algname{MbSGD}. To the best of our knowledge, this is the first speedup analysis of \algname{LocalSGD} for non-convex functions. Moreover, this result should not be considered a trivial extension of the existing convex speedup analysis (cf. \cref{lemma:LocalSGD faster convex previous}) because we do not rely on uniform gradient similarity here. 
\end{remark} 

As a small detour, we would like to briefly present a convex result here (although the focus of this work is non-convex analysis). We show that, for convex functions, following similar techniques used in our \cref{thm:LocalSGD faster}, we can strictly improve the previous results of \cref{lemma:LocalSGD faster convex previous}. 
    \begin{theorem}
    \label{thm:LocalSGD faster convex}
    Under \cref{assumption:GS}, if all the local functions $f_i$ are convex, $\xx^* \in \argmin _{\xx \in \R^d} {f(\xx)}$, and there exists some $D\ge 0$ such that $\norm{ \bar \xx_0 - \xx^*} \le D$, then there exists $\eta > 0$ such that \algname{LocalSGD} ensures the following upper bound on $\frac 1 T \sum _{t=0} ^{T-1} \Eb { f \left( \bar \xx_t \right) } - f^*$: 
    \begin{equation}
    \label{eq:LocalSGD convex faster rate}
    \cO \left( \frac {LD^2} {{ \tau} R} + \frac {\sigma D} {\sqrt{n \tau R}} + \left( \frac {L { \zeta}^2 D^4} {R^2} \right) ^{\frac 1 3} + \left( \frac {L \sigma^2 D^4} {\tau R^2} \right) ^{\frac 1 3} \right) . 
    \end{equation} 
    \end{theorem}
    \begin{remark}
    As compared to the previous speedup result \eqref{eq:LocalSGD convex faster rate previous}, the heterogeneity term in \eqref{eq:LocalSGD convex faster rate} replaces the dependence on uniform gradient similarity with \emph{standard} gradient similarity, which is always tighter, and can be substantially tighter up to a factor of $n^{\frac 1 3}$ with imbalanced heterogeneity. 
    Moreover, \cref{thm:LocalSGD faster convex} shows that standard gradient similarity ${ \zeta}^2 = \cO (1/R)$ is sufficient to achieve speedup for \algname{LocalSGD}. In contrast, the stronger condition of uniform gradient similarity \( { \barzeta}^2 = \cO \left( 1/R \right) \) in~\cite{woodworth2020minibatch,patel2024limits} is \emph{redundant}. 
    \end{remark}

\noindent{\bf Novelty.}
    The key technique in our proof of \cref{thm:LocalSGD faster,thm:LocalSGD faster convex} is a ``variance trick''. A crucial term in the analysis is the variance between the workers, \ie, 
    \[
    \frac 1 n \sum _{i=1} ^n \norm {\xx_t^i - \bar \xx_t}^2 := \Xi_t, \text{ where } t=r\tau + k, \ k \in [\tau-1] . 
    \]
    Essentially, many prior works either upper bound $\Xi_t$ by $\frac 1 n \sum _{i=1} ^n \norm {\xx_t^i - \bar \xx_{r\tau}}^2$~\citep{koloskova2020unified,karimireddy2020scaffold}, which ignores the worker similarity; or upper bound $\Xi_t$ by $\sup _{i,j} \norm {\xx_t^i - \xx_t^j}^2$~{\citep{woodworth2020minibatch,patel2024limits}}, which loses the symmetry. We instead upper bound $\Xi_t$ by $\frac 1 n \sum _{i=1} ^n \norm {\xx_t^i - \bar \xx_{t-1} + \eta \nabla f(\bar \xx_{t-1})}^2$ and then unroll the recursion. This small trick simplifies the analysis and leads to a strict improvement over previous analyses. 
\qed

\begin{theorem}
\label{thm:LocalSGD improved conditioning}
    Under \cref{assumption:GS,assumption:HS,assumption:Lipschitz Hessian}, there exists $\eta>0$ such that \algname{LocalSGD} ensures the following upper bound on $\frac 1 T \sum _{t=0}^{T-1} \E \norm {\nabla f(\bar \xx_t)}^2$: 
    \begin{equation} 
    \label{eq:LocalSGD HS rate}
    \begin{aligned}
    \cO \Bigg( \frac {L\Delta} {R} & + \sqrt {\frac {L\Delta \sigma^2} {n\tau R}} + \left( \frac {{ \bardelta}\Delta \zeta} {R} \right)^{\frac23} \\ 
    & \quad + \frac {(L\Delta\sigma)^{\frac23}} {\tau^{\frac13}R^{\frac23}} + \left( \frac {\cM^2\Delta^4\zeta^4} {R^4} \right)^{\frac15}  \Bigg) . 
    \end{aligned}
    \end{equation}
\end{theorem} 

\begin{remark}
\label{remark:LocalSGD HS}
As compared to the previous non-convex analysis (\cref{lemma:classic LocalSGD}), \cref{thm:LocalSGD improved conditioning} shows that, under the assumptions of uniform Hessian similarity and Lipschitz continuous Hessian, it is possible to improve the conditioning in the heterogeneity term. 
As a simplification of \eqref{eq:LocalSGD HS rate}, for the case $L=\Delta=\sigma=1$ and $\tau=\infty$, the asymptotic rate of \algname{LocalSGD} is 
\[
\cO \left( \frac 1 R + \left( \frac {{ \bardelta}\zeta} {R} \right) ^{2/3} + \left( \frac {\cM^2\zeta^4} {R^4} \right)^{1/5} \right) . 
\]

We see that the above rate of \algname{LocalSGD} is still no faster than that of \algname{MbSGD}, but the conditioning required to match \algname{MbSGD} has been relaxed. 
Asymptotically, this rate of \eqref{eq:LocalSGD HS rate} can match that of \algname{MbSGD} under the conditioning of 
\(
    { \bardelta^2} \zeta^2 + \cM^2 \zeta^4 = \cO \left( 1 / R \right) . 
\)
For $\cM^2 = \cO \left( R \right)$, this is weaker than the conditioning of $\zeta^2 = \cO (1 / R)$ in the classic analysis (\cf \cref{remark:classic LocalSGD}). 

We also see that when $\zeta^2 = \varOmega \left( 1/R \right)$, the heterogeneity term dominates the asymptotic rate in the classic analysis (\cf \cref{remark:classic LocalSGD}). Meanwhile, if $\cM^3 \zeta = \cO \left( R \right)$, with the new analysis, our rate in \eqref{eq:LocalSGD HS rate} demonstrates a strict improvement over the previous analysis. This is a perhaps surprising result showing that the local gradient steps can benefit from higher-order conditions, even without variance reduction. 
\end{remark}

\noindent{\bf Novelty.}
    We outline the proof techniques for \cref{thm:LocalSGD improved conditioning}. At each iteration $t$, the desired step is $\nabla f(\bar \xx_t)$, but the actual step is approximately $\frac 1 n \sum _{i=1} ^n \nabla f_i (\xx_t^i)$. Thus, a crucial part of the analysis is bounding the discrepancy between these two, \ie, 
    \[
    \norm {\frac 1 n \sum _{i=1} ^n \nabla f_i (\xx_t^i) - \nabla f(\bar \xx_t)}^2 := Q_t .
    \]
    Previous approaches, such as~\cite{koloskova2020unified,karimireddy2020scaffold}, simply upper bound $Q_t$ by $L^2 \Xi_t$. Instead, under \cref{assumption:HS,assumption:Lipschitz Hessian}, we derive a tighter upper bound of \(8 { \bardelta}^2 \Xi_t + \frac {\cM^2} 2 {\Xi_t}^2 \) on $Q_t$.
    This technique alone does not suffice due to the stochastic noise in ${{\Xi_t}^2}$. To this end, we combine it with the construction of a ``noiseless sequence'', which helps to decouple the stochasticity from the gradient discrepancy, and the desired result finally follows. 
\qed

\subsection{{New Analysis of SCAFFOLD}}
\label{sec:scaffold}

Currently, the only speedup analysis of \algname{SCAFFOLD} requires all $f_i$ functions to be quadratic. In this section, we show the speedup of \algname{SCAFFOLD}, while relaxing the ``all-quadratic'' restriction. 

\begin{theorem}
\label{thm:SCAFFOLD speedup}
    Under \cref{assumption:mean HS,assumption:weak convexity}, there exists $\eta>0$ such that \algname{SCAFFOLD} ensures the following upper bound on $\frac 2 T \sum _{r=0}^{R-1} \sum _{k=0}^{\tau-1} \E \norm {\nabla f(\bar \xx_{2r \tau + \tau + k})}^2$: 
    \begin{equation}
    \label{eq:improved general SCAFFOLD rate}
    \cO \left( \left( { \frac L \tau + \sqrt{L\delta} + \rho} \right) \frac {\Delta} {R} + \sqrt { \frac {L \Delta \sigma^2} {n \tau R} } + \frac {(L\Delta\sigma)^{\frac 2 3}} {\tau^{\frac 1 3}R^{\frac 2 3}}  \right) . 
    \end{equation}
\end{theorem} 

\begin{remark}
As compared to the rate of \algname{MbSGD}~(\cref{lemma:MbSGD}) and the previous rate of \algname{SCAFFOLD}~(\cref{lemma:original SCAFFOLD}), our new rate of \algname{SCAFFOLD} in \cref{thm:SCAFFOLD speedup} has a faster optimization term of $\cO \left( (\frac L \tau + \sqrt {L \delta} + \rho) \frac \Delta R \right)$. 
As a simplification of \eqref{eq:improved general SCAFFOLD rate}, for the case $L=\Delta=\sigma=1$ and $\tau=\infty$, the asymptotic rate of \algname{SCAFFOLD} is 
\(
\cO \left( \frac {{ \sqrt{\delta} + \rho}} R \right) , 
\)
which is faster than that of \algname{MbSGD}. 
\end{remark}

\noindent{\bf Novelty.} 
    This result extends the speedup analysis of \algname{SCAFFOLD} beyond quadratic cases, addressing a challenging open problem stated in~\cite{karimireddy2020scaffold}. We also weaken the uniform Hessian similarity required in prior work to the standard one.  
\qed

Finally, for theoretical purposes, we show that, under \cref{assumption:Lipschitz Hessian} with $\cM=0$, 
we can get more benefit from Hessian similarity, thus improving the optimization term from $\cO \left( (\frac L \tau + \sqrt {L \delta} + \rho) \frac \Delta R \right)$ to $\cO \left( (\frac L \tau + \sqrt {\bardelta \delta} + \rho) \frac \Delta R \right)$.  
\begin{theorem}
\label{thm:Lipschitz Hessian SCAFFOLD}
    Under Assumptions~\ref{assumption:mean HS},~\ref{assumption:HS},~\ref{assumption:weak convexity}~and~\ref{assumption:Lipschitz Hessian} with $\cM=0$, there exists $\eta > 0$ \st \algname{SCAFFOLD} ensures the following upper bound on $\frac 2 T \sum _{r=0}^{R-1} \sum _{k=0}^{\tau-1} \E \norm {\nabla f(\bar \xx_{2r \tau + \tau + k})}^2$: 
    \[
    \cO \left( \left( { \frac L \tau + \sqrt {\bardelta \delta} + \rho} \right) \frac {\Delta} {R} + \sqrt { \frac {L \Delta \sigma^2} {n \tau R} } + \frac {({ \bardelta} \Delta\sigma)^{\frac 2 3}} {\tau^{\frac 1 3}R^{\frac 2 3}}  \right) . 
    \]
\end{theorem}

\section{Synthetic Experiments}
\label{sec:expr}

\noindent{\bf Experiment setup.}
We consider a distributed setting with $n$ workers, where the $i$-th worker holds a dataset $(\mA_i, \yy_i) \in \R^{m_i \times d} \times \R^{m_i}$, for each $i \in [n]$. The total number of data points across all workers is denoted by $m = \sum _{i=1} ^n m_i$. 
We focus on solving the following regression problem with smoothed Huber loss and non-convex regularizer: 
\begin{equation}
\label{eq:logistic_regression}
    \min _{\xx \in \R^d} f(\xx) = \frac 1 n \sum_{i=1} ^n \left[ \cL (\xx; \mA_i, \yy_i) + \lambda \cdot \sum _{l=1} ^d \frac {\xx_l^2} {1 + \xx_l^2} \right] , 
\end{equation}
where 
\(
\cL(\xx; \mA_i, \yy_i) = \frac n {m} \cdot \sum _{j=1} ^{m_i} h \left( 
\mA_i(j) \cdot \xx - \yy_i(j) \right) ,  
\) 
in which $h(\cdot)$ is defined as follows: 
\[
h (u) = \begin{cases} 
    \frac 1 2 u^2, & |u| \le 1, \\
    - \frac 1 6 (|u|-1)^3 + \frac 1 2 u^2, & 1 < |u| \le 2, \\
    \frac 3 2 |u| - \frac 7 6, & |u| > 2. 
\end{cases}
\]
In all experiments, we set the dimension $d=100$ and the number of workers $n=10$. We fix the regularization parameter $\lambda=0.01$ and set the following conditionings: $L\approx 1$, $\Delta \approx 1$ and $\sigma \approx 0.01$. We tune the stepsizes over $\{0.0003, 0.001, 0.003, 0.01, 0.03, 0.1, 0.3\}$.

In this section, we present the experiments that validate the dependence on Hessian similarity $\delta$. Further details and additional experiments on the other parameters can be found in \cref{sec:expr details}. 

\begin{figure}[ht]
    \centering
    \subfigure{
        \includegraphics[width=0.22\textwidth]{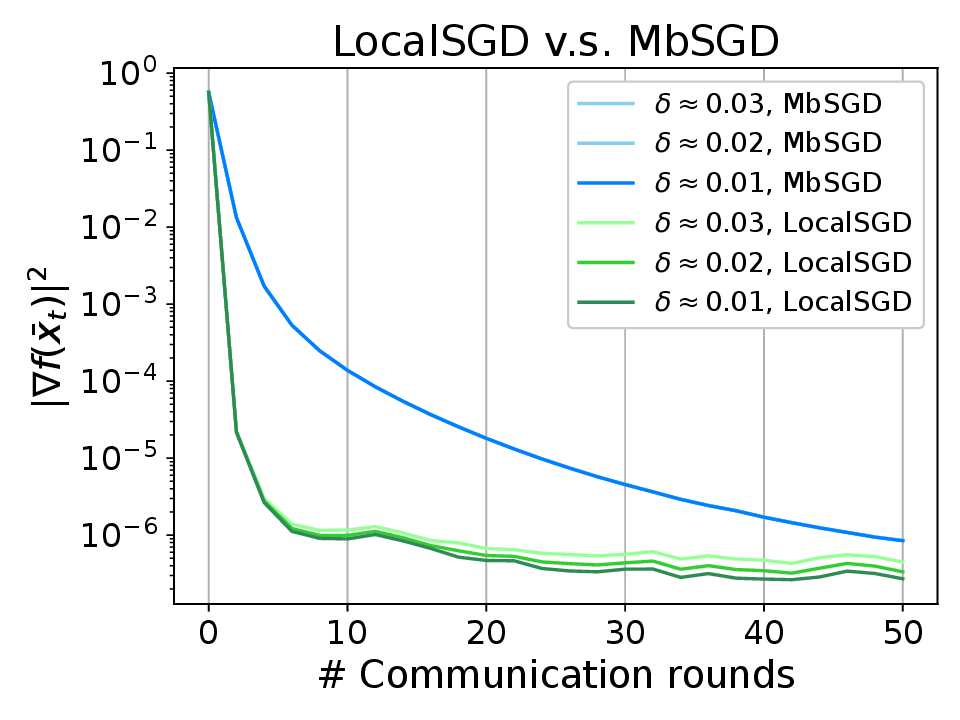} 
    }
    \subfigure{
        \includegraphics[width=0.22\textwidth]{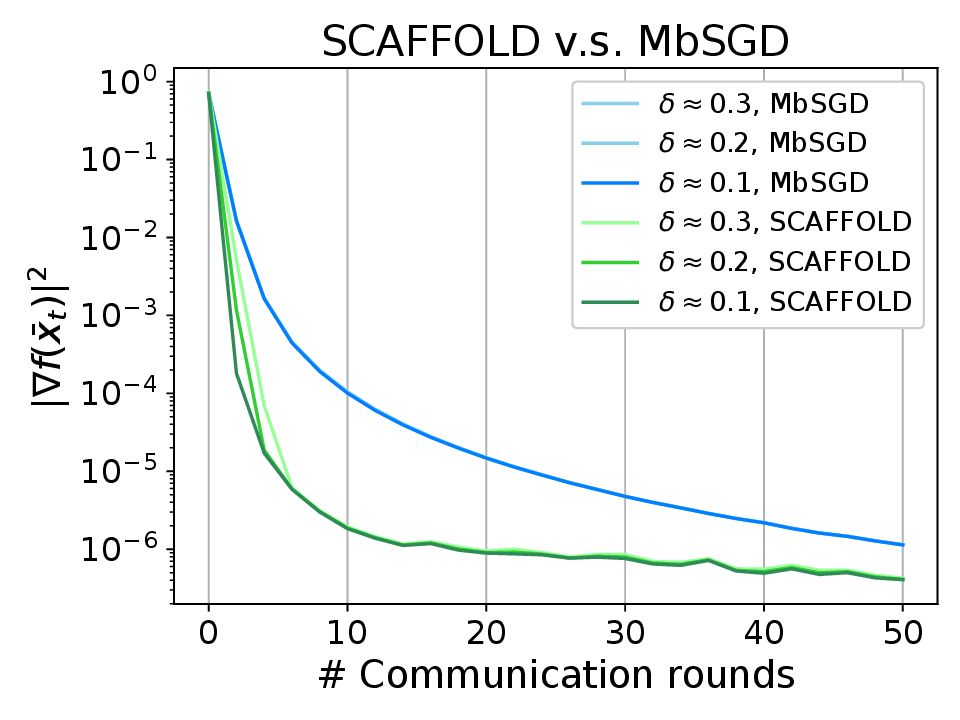}  
    }
    \vspace{-1em}
    \caption{Comparisons between the convergence rates.}
    \label{fig:localsgd/scaffold_vs_mbsgd}
\end{figure}

\textbf{\algname{LocalSGD} v.s. \algname{MbSGD}.} The theoretical results in~\cref{thm:LocalSGD faster,thm:LocalSGD improved conditioning} imply that \algname{LocalSGD} can achieve a speedup over \algname{MbSGD} for weakly convex functions and benefit from higher order conditions. To validate these results, we generated data such that $\zeta \approx 0.03$ and  $\delta \approx 0.01, 0.02, 0.03$, and set the communication interval $\tau=50$. As shown in~\cref{fig:localsgd/scaffold_vs_mbsgd}~(left), \algname{MbSGD} converges at the same rate despite variations in $\delta$. In contrast, \algname{LocalSGD} converges faster as $\delta$ decreases. The most pronounced differences between the \algname{LocalSGD} runs are observed in later communication rounds ($\ge 10$), where the heterogeneity term dominates the rate.

\textbf{\algname{SCAFFOLD} v.s. \algname{MbSGD}.} \Cref{thm:SCAFFOLD speedup} shows the speedup of \algname{SCAFFOLD} from Hessian similarity and weak convexity. To validate this result, we generated data such that $\delta \approx 0.1, 0.2, 0.3$, and set the communication interval $\tau=50$. As shown in~\cref{fig:localsgd/scaffold_vs_mbsgd}~(right), \algname{MbSGD} converges at the same rate despite variations in $\delta$. In contrast, \algname{SCAFFOLD} converges faster as $\delta$ decreases. The most pronounced difference between the \algname{SCAFFOLD} runs is observed in early communication rounds ($\le 10$), where the optimization term dominates the rate.

\section{Limitations and Future Work}

This paper revisits the convergence rates of \algname{LocalSGD} and \algname{SCAFFOLD} under gradient similarity, Hessian similarity, weak convexity, and Lipschitz Hessian. Our extensive studies demonstrate improved rates for both algorithms across various settings, suggesting that our new analysis can better capture their effectiveness as observed in real-world applications.

While our analysis focuses on the standard algorithms in the intermittent communication setting, many orthogonal techniques have been applied to the framework to further reduce communication cost. These include the use of previous updates for variance reduction~\citep{karimireddy2020scaffold}, quantization and sparsification of the updates~\citep{alistarh2017qsgd,stich2018sparsified}, decentralized communication network~\citep{lian2017can,sadiev2022decentralized}, partial participation~\citep{mcmahan2017communication}, and asynchronous communication~\citep{recht2011hogwild,nguyen2018sgd}. We do not cover these techniques in our analysis and leave them as future work. 


\subsubsection*{Acknowledgements}
The authors thank for the helpful discussions with Eduard Gorbunov, Kumar Kshitij Patel, Anton Rodomanov, and Ali Zindari during the preparation of this work. 

This work was partially done during the first author's stays at CISPA and at MBZUAI. The first author also acknowledges ERC CoG 863818 (ForM-SMArt) and Austrian Science Fund (FWF) 10.55776/COE12. 

\bibliography{reference}

 \onecolumn
 \appendix

 \hsize\textwidth
  \linewidth\hsize \toptitlebar {\centering
  {\Large\bfseries Supplementary Materials for Revisiting LocalSGD and SCAFFOLD \par}}
 \bottomtitlebar 
 
 \section{Proof Details}
 \label{appendix_sec:proof_details}
 
 \subsection{Technical Lemmas} 
 
 In this section, we present several technical lemmas and remark on their usages in our analysis. 
 
 For simplicity, we denote $\bar \xx = \frac 1 n \sum _{i=1} ^n \xx^i$ and $\Xi = \frac 1 n \sum_{i=1} ^n \norm {\xx^i - \bar\xx}^2$. We assume all the vectors are in $\R^d$. 
 
 \begin{lemma}[variance trick]
     For any $\yy \in \R^d$, 
     \begin{equation}
     \label{eq:variance relax}
     \Xi \le \sum_{i=1} ^n \norm {\xx^i - \yy}^2 .
     \end{equation}
 \end{lemma}
 \begin{proof}
     Since $\sum_{i=1} ^n \norm {\xx^i - \yy}^2$ is strongly convex in $\yy$, let $\yy^* = \arg \min _{\yy \in \R^d} \sum_{i=1} ^n \norm {\xx^i - \yy}^2$. Then,
     \[
     \left. \frac {d \left( \sum_{i=1} ^n \norm {\xx^i - \yy}^2 \right)} {d \yy} \right | _{\yy^*} = \left. 2 (\yy - \bar \xx) \right | _{\yy^*} = \0, 
     \]
     and we have $\yy^* = \bar \xx$. 
 \end{proof}
 \begin{remark}
     We will use this variance trick frequently in our analysis. For instance, when upper bounding $\Xi_{t+1}$ recursively, we can use the following inequality:
     \[
     \Xi_{t+1} \le \frac 1 n \sum _{i=1} ^n \norm {\xx_{t+1}^i - \bar \xx_t + \eta \nabla f(\bar \xx_t)} ^2 . 
     \]
 \end{remark}

 \begin{lemma}
 \label{lemma:refined upper bound of Q}
     Under \cref{assumption:HS,assumption:Lipschitz Hessian}, we have
     \begin{equation}
     \label{eq:refined upper bound of Q}
     \norm { \frac 1 n \sum_{i=1} ^n \nabla f_i(\xx^i) - \nabla f(\bar \xx) }^2 \le 8\bardelta^2\Xi + \frac {\cM^2} 2 \Xi^2 .
     \end{equation}
 \end{lemma}
 \begin{proof} 
     Let $\hat f \in \mathbf{conv} \{ f_1, \cdots, f_n \}$ \st $\hat f$ has $\cM$-Lipschitz continuous Hessian. 
     We have 
     \begin{equation}
     \label{eq:norm of average gradient difference}
     \begin{aligned}
     &\quad \norm { \frac 1 n \sum_{i=1} ^n \nabla f_i(\xx^i) - \nabla f(\bar \xx) }^2 \\
     &= \norm { \frac 1 n \sum_{i=1} ^n \left( \nabla f_i(\xx^i) - \nabla f_i(\bar \xx) \right) }^2 \\
     &\le 2 \Bigg( \norm { \frac 1 n \sum_{i=1} ^n \nabla \hat f(\xx^i) - \nabla \hat f(\bar \xx) }^2 
     + \norm { \frac 1 n \sum_{i=1} ^n \left( \nabla f_i(\xx^i) - \nabla \hat f(\xx^i) - \nabla f_i(\bar \xx) + \nabla \hat f(\bar \xx) \right) }^2 \Bigg) . 
     \end{aligned} 
     \end{equation}
 
     We then upper bound the above two terms respectively: 
     \begin{equation}
     \label{eq:upper bound average gradient difference on reference function}
     \begin{aligned}
     &\quad \norm { \frac 1 n \sum_{i=1} ^n \nabla \hat f(\xx^i) - \nabla \hat f(\bar \xx) }^2 \\
     &= \norm {\frac 1 n \sum_{i=1} ^n \int _0^1 \nabla^2 \hat f \left( \bar \xx + u(\xx^{i} - \bar \xx) \right) (\xx^{i} - \bar \xx) du }^2 \\
     &= \norm {\frac 1 n \sum_{i=1} ^n \int _0^1 \left[ \nabla^2 \hat f \left( \bar \xx + u(\xx^{i} - \bar \xx) \right) - \nabla^2 \hat f(\bar \xx) \right] (\xx^{i} - \bar \xx) du}^2 \\
     &\le \left[ \frac 1 n \sum_{i=1} ^n \int _0^1 \norm { \nabla^2 \hat f \left( \bar \xx + u(\xx^{i} - \bar \xx) \right) - \nabla^2 \hat f(\bar \xx) } \cdot \norm { \xx^{i} - \bar \xx } du \right] ^2 \\
     &\refLE{eq:Lipschitz Hessian} \left[ \frac 1 n \sum_{i=1} ^n \int _0^1 \cM u \norm{\xx^{i} - \bar \xx} ^2 du \right]^2 \\
     &= \frac {\cM^2} 4 \Xi^2 ,
     \end{aligned}
     \end{equation}
     and
     \begin{equation}
     \label{eq:upper bound hessian similarity}
     \begin{aligned}
     &\quad \norm { \frac 1 n \sum_{i=1} ^n \left( \nabla f_i(\xx^i) - \nabla \hat f(\xx^i) - \nabla f_i(\bar \xx) + \nabla \hat f(\bar \xx) \right) }^2 \\
     &\le \frac 1 n \sum _{i=1} ^n \norm { \nabla f_i(\xx^i) - \nabla \hat f(\xx^i) - \nabla f_i(\bar \xx) + \nabla \hat f(\bar \xx) } ^2 \\
     &{\le} \frac 1 n \sum _{i=1} ^n \left( 2 \bardelta \norm { \xx^i - \bar \xx } \right) ^2 \\
     &= 4 \bardelta^2 \Xi . 
     \end{aligned}
     \end{equation}
 
     Finally, \eqref{eq:refined upper bound of Q} follows from plugging \eqref{eq:upper bound average gradient difference on reference function} and \eqref{eq:upper bound hessian similarity} into \eqref{eq:norm of average gradient difference}. 
 \end{proof} 
 \begin{remark}
     Lemma~\ref{lemma:refined upper bound of Q} is the key to the benefit from higher order conditions in Theorem~\ref{thm:LocalSGD improved conditioning} and Theorem~\ref{thm:Lipschitz Hessian SCAFFOLD}. Instead, without higher order conditions, we can only upper bound the term as follows: 
     \begin{equation}
     \label{eq:normal upper bound of Q}
     \norm { \frac 1 n \sum_{i=1} ^n \nabla f_i(\xx^i) - \nabla f(\bar \xx) }^2 \le \frac 1 n \sum_{i=1} ^n \norm { \nabla f_i(\xx^i) - \nabla f_i(\bar \xx) }^2 \le L^2\Xi .
     \end{equation}
 \end{remark}
 
 \begin{lemma}
     Under \cref{assumption:weak convexity}, we have 
     \begin{equation}
     \label{eq:smooth weakly convex inequality}
     (L - \rho) \lin { \xx - \yy, \nabla f_i(\xx) - \nabla f_i(\yy) } \ge \norm {\nabla f_i(\xx) - \nabla f_i(\yy)}^2 - \rho L \norm {\xx - \yy}^2 .
     \end{equation}
 \end{lemma}
 \begin{proof}
     $f_i(\xx) + \frac {\rho} 2 \xx \T \xx$ is $(L+\rho)$-smooth and convex. By $\textnormal{(2.1.8)}$ in~\cite{nesterov2003introductory}, we have
     \[
     \frac {1} {L+\rho} \norm { \nabla f_i(\xx) + \rho \xx - \nabla f_i(\yy) - \rho \yy }^2 \le \lin { \nabla f_i(\xx) + \rho \xx - \nabla f_i(\yy) - \rho \yy , \xx - \yy } .
     \]
     Then, after multiplying both sides by $(L+\rho)$, a rearrangement will yield \eqref{eq:smooth weakly convex inequality}. 
     
 \end{proof}
 
 \begin{lemma}
 \label{lemma:weak convexity contraction}
     Under \cref{assumption:weak convexity}, for $\eta \le \frac 2 {L-\rho}$, we have 
     \begin{equation}
     \label{eq:weak convexity contraction}
     \norm { (\xx - \eta \nabla f_i(\xx)) - (\yy - \eta \nabla f_i(\yy)) }^2 \le \left( 1 + \frac {2L} {L-\rho} \cdot \rho \eta \right) \norm {\xx - \yy}^2 .
     \end{equation}
 \end{lemma}
 \begin{proof}
     For $\eta \le \frac 2 {L-\rho}$, 
     \[
     \begin{aligned}
     &\quad \norm { \xx - \eta \nabla f_i(\xx) - (\yy - \eta \nabla f_i(\yy)) }^2 \\
     &= \norm {\xx - \yy}^2 - 2\eta \lin { \xx - \yy, \nabla f_i(\xx) - \nabla f_i(\yy) } + \eta^2 \norm { \nabla f_i(\xx) - \nabla f_i(\yy) } ^2 \\
     &\refLE {eq:smooth weakly convex inequality} \left( 1 + \frac {2L} {L-\rho} \cdot \rho \eta \right) \norm {\xx - \yy}^2 - \left( \frac {2\eta} {L - \rho} - \eta^2 \right) \norm { \nabla f_i(\xx) - \nabla f_i(\yy) }^2 \\
     &\le \left( 1 + \frac {2L} {L-\rho} \cdot \rho \eta \right) \norm {\xx - \yy}^2 .
     \end{aligned}
     \]
 \end{proof}
 \begin{remark}
     Lemma~\ref{lemma:weak convexity contraction} is the key to prove the speedup from weak convexity. 
     Instead, without weak convexity, we can only upper bound the term as follows: 
     \begin{equation}
     \label{eq:smooth contraction}
     \begin{aligned}
     \norm { \xx - \eta \nabla f_i(\xx) - (\yy - \eta \nabla f_i(\yy)) }^2 
     &\le \left( \norm {\xx - \yy} + \eta \norm {\nabla f_i(\xx) - \nabla f_i(\yy)} \right) ^2 \\
     &\le (1+L\eta)^2 \norm {\xx - \yy}^2 .  
     \end{aligned}
     \end{equation}
     This upper bound in \cref{eq:smooth contraction} is looser compared to the one in \cref{eq:weak convexity contraction}. 
 \end{remark}
 
 \subsection{The Analysis of LocalSGD} 
 \label{sec:LocalSGD proofs}
 
 We recall the following notations: for $t=r\tau+k$, $r\in [0,R-1]$, $k\in [0,\tau-1]$, we denote 
 \[
 r(t) = r\tau, \quad \Xi_t = \frac 1 n \sum_{i=1} ^n \norm {\xx_t^i - \bar\xx_t}^2, \quad Q_t = \norm { \frac 1 n \sum_{i=1} ^n \nabla f_i(\xx_t^i) - \nabla f(\bar \xx_t) }^2 .
 \] 
 
 Define the following \emph{noiseless sequence}: for $t \in [0, T-1]$, 
 \[
 \hat \xx _{t} ^i = \begin{cases}
     \bar \xx_{t} , & \text{ if } t \text{ is a multiple of } \tau, \\ 
     \hat \xx_{t-1}^i - \eta \nabla f_i (\hat \xx_{t-1}^i), & \text{ otherwise} ,
     \end{cases} 
 \]
 and 
 \[
 \hat \xx_t = \frac 1 n \sum _{i=1} ^n \hat \xx_t^i , \quad \hat Q_t = \frac 1 n \sum _{i=1} ^n \norm { \nabla f_i (\hat \xx_t^i) - \nabla f (\hat \xx_t) }^2 , \quad \hat \Xi_t = \frac 1 n \sum_{i=1} ^n \norm { \hat \xx_t^i - \hat \xx_t }^2 .
 \]
 
 \begin{lemma}
   For $\eta \le \frac 1 L$, \algname{LocalSGD} ensures
   \begin{equation}
   \label{eq:localsgd descent lemma}
     \frac 1 T \sum _{t=0} ^{T-1} \E \norm {\nabla f(\bar \xx_t)}^2 \le \frac {2 \Delta} {\eta T} + \eta \frac {L \sigma^2} {n} + \frac {1} {T} \sum _{r=0} ^{R-1} \sum _{k=0} ^{\tau-1} \left( {6 L^2 \frac 1 n \sum _{i=1} ^n \Eb { \norm { \xx_{r \tau + k}^i - \hat \xx_{r \tau + k}^i }^2 } + 3 \Eb {\hat Q _{r \tau + k}}} \right) . 
   \end{equation}
 \end{lemma}
 \begin{proof} 
   For $\eta \le \frac 1 L$, 
   \[
   \begin{aligned}
   &\quad \Eb {f(\bar \xx_{t+1})} \\
   &\refLE{eq:borel} \Eb {f(\bar \xx_t)} - \eta \E \lin { \nabla f(\bar \xx_t) , \frac 1 n \sum _{i=1} ^n \nabla f_i(\xx_t^i) } + \eta^2 \frac L 2 \E \norm { \frac 1 n \sum _{i=1} ^n \nabla f_i(\xx_t^i) }^2 + \eta^2 \frac {L\sigma^2} {2n} \\
   &= \Eb {f(\bar \xx_t)} - \frac {\eta} 2 \E \norm {\nabla f(\bar \xx_t)}^2 - \left( \frac \eta 2 - \eta^2 \frac L 2 \right) \E \norm { \frac 1 n \sum _{i=1} ^n \nabla f_i(\xx_t^i) }^2 + \eta^2 \frac {L\sigma^2} {2n} + \frac {\eta} 2 \Eb { Q_t } \\
   &\le \Eb {f(\bar \xx_t)} - \frac {\eta} 2 \E \norm {\nabla f(\bar \xx_t)}^2 + \eta^2 \frac {L\sigma^2} {2n} + \frac {\eta} 2 \Eb { Q_t } . 
   \end{aligned}
   \] 
   Summing over $t$ from $0$ to $T-1$ and multiplying both sides by $\frac 2 {\eta T}$, we have 
     \[ 
     \begin{aligned}
     \frac 1 T \sum _{t=0} ^{T-1} \E \norm {\nabla f(\bar \xx_t)}^2 &\le \frac {2} {\eta T} \E \left[ f(\bar \xx_0) - f(\bar \xx_T) \right] + \eta \frac {L \sigma^2} {n} + \frac {1} {T} \sum _{r=0} ^{R-1} \sum _{k=0} ^{\tau-1} \Eb {Q_{r\tau+k}} \\ 
     &\le \frac {2 \Delta} {\eta T} + \eta \frac {L \sigma^2} {n} + \frac {1} {T} \sum _{r=0} ^{R-1} \sum _{k=0} ^{\tau-1} \Eb {Q_{r\tau+k}} . 
     \end{aligned}
     \]
   Finally, \cref{{eq:localsgd descent lemma}} follows from the following inequality: 
   \[
     \begin{aligned}
         \Eb {Q_t} &\le 3 \Eb { \norm { \frac 1 n \sum _{i=1} ^n \left[ \nabla f_i (\xx_t^i) - \nabla f_i (\hat \xx_t^i) \right] }^2 } + 3 \Eb { \norm { \nabla f (\bar \xx_t) - \nabla f (\hat \xx_t) }^2 } + 3 \Eb { \hat Q_t } \\
         &\le \frac 3 n \sum _{i=1} ^n \Eb { \norm { \nabla f_i (\xx_t^i) - \nabla f_i (\hat \xx_t^i) }^2 } + 3 \Eb { \norm { \nabla f (\bar \xx_t) - \nabla f (\hat \xx_t) }^2 } + 3 \Eb { \hat Q_t } \\
         &\le 3 L^2 \frac 1 n \sum _{i=1} ^n \Eb { \norm { \xx_t^i - \hat \xx_t^i }^2 } + 3 L^2 \norm { \bar \xx_t - \hat \xx_t }^2 + 3 \Eb {\hat Q_t} \\
         &\le 6 L^2 \frac 1 n \sum _{i=1} ^n \Eb { \norm { \xx_t^i - \hat \xx_t^i }^2 } + 3 \Eb {\hat Q_t} .
     \end{aligned}
   \]
 
 \end{proof}
 
 \begin{lemma}[variance trick]
     For any $\gamma > 0$, when $(t+1)$ is not a multiple of $\tau$, we have 
     \begin{equation}
     \label{eq:variance trick in localsgd}
     {\hat \Xi_{t+1}} \le \frac 1 n \sum _{i=1} ^n (1 + \gamma) \norm { \hat \xx_{t}^i - \eta \nabla f_i(\hat \xx_t^i) - \left( \hat \xx_t - \eta \nabla f_i(\hat \xx_t) \right) } ^2 + (1+\gamma^{-1}) \eta^2 \zeta^2 . 
     \end{equation}
 \end{lemma}
 \begin{proof}
     For any $\gamma > 0$, 
     \begin{equation}
     \label{eq:variance trick proof in localsgd}
     \begin{aligned}
         &\quad {\hat \Xi_{t+1}} \\
         &\refLE{eq:variance relax} \frac 1 n \sum _{i=1} ^n \norm { \hat \xx_{t+1}^i - \left( \hat \xx_t - \eta \nabla f(\hat \xx_t) \right) } ^2 \\
         &= \frac 1 n \sum _{i=1} ^n \norm { \hat \xx_{t}^i - \eta \nabla f_i(\hat \xx_t^i) - \left( \hat \xx_t - \eta \nabla f(\hat \xx_t) \right) } ^2 \\
         &\refLE{eq:GS bound} \frac 1 n \sum _{i=1} ^n (1 + \gamma) \norm { \hat \xx_{t}^i - \eta \nabla f_i(\hat \xx_t^i) - \left( \hat \xx_t - \eta \nabla f_i(\hat \xx_t) \right) } ^2 + (1+\gamma^{-1}) \eta^2 \zeta^2 .  
     \end{aligned}
     \end{equation}
 \end{proof}
 
 \subsubsection{Proof of Theorem~\ref{thm:LocalSGD faster}}
 \label{sec:LocalSGD faster proof}
 
 \begin{lemma}
     Under \cref{assumption:weak convexity}, 
     for 
     \(
     \eta \le \min \left\{ \frac 2 {L-\rho} , \frac {1 - {\rho} / L} {4 \rho (\tau-1)} \right\} , 
     \)
     \algname{LocalSGD} ensures 
     \begin{equation}
     \label{eq:normal distance lemma LocalSGD distance to}
     \E \norm { \xx_t^i - \hat \xx_t^i }^2 \le 3 k \eta^2 \sigma^2 , 
     \end{equation}
     where $k = t - r(t)$. 
 \end{lemma} 
 The proof of the above lemma is similar to the proof of \cref{lemma:localsgd variance lemma}. \cref{lemma:weak convexity contraction} is used in this proof when bounding the distance recursively. 
 
 \begin{lemma}
 \label{lemma:weak convexity distance lemma LocalSGD noiseless}
     Under \cref{assumption:GS,assumption:weak convexity},  
     for 
     \(
     \eta \le \min \left\{ \frac 2 {L-\rho} , \frac {1 - {\rho} / L} {6 \rho (\tau-1)} \right\} , 
     \)
     \algname{LocalSGD} ensures 
     \begin{equation}
     \label{eq:normal distance lemma LocalSGD}
     \hat \Xi_t \le 9 k (\tau-1) \eta^2 \zeta^2 , 
     \end{equation}
     where $k = t - r(t)$.
 \end{lemma}
 The proof of the above lemma is similar to the proof of \cref{lemma:normal distance lemma LocalSGD noiseless}. \cref{lemma:weak convexity contraction} is used in this proof when bounding the distance recursively. 
 
 Now we are ready to prove Theorem~\ref{thm:LocalSGD faster}. 
 \begin{proof}[Proof of Theorem~\ref{thm:LocalSGD faster}] 
     For $\eta \le \min \left\{ \frac 1 L, \frac {1 - {\rho} / L} {6 \rho (\tau-1)} \right\}$, 
     \[
     \begin{aligned}
     &\quad \frac 1 T \sum _{t=0} ^{T-1} \E \norm {\nabla f(\bar \xx_t)}^2 \\ 
     &\refLE{eq:localsgd descent lemma} \frac {2 \Delta} {\eta T} + \eta \frac {L \sigma^2} {n} + \frac {1} {T} \sum _{r=0} ^{R-1} \sum _{k=0} ^{\tau-1} \left( {6 L^2 \frac 1 n \sum _{i=1} ^n \Eb { \norm { \xx_{r \tau + k}^i - \hat \xx_{r \tau + k}^i }^2 } + 3 \Eb {\hat Q _{r \tau + k}}} \right) \\
     &\stackrel {\eqref{eq:normal distance lemma LocalSGD distance to} \eqref{eq:normal upper bound of Q}} {\le} \frac {2 \Delta} {\eta T} + \eta \frac {L \sigma^2} {n} + { 9 (\tau-1) L^2 \eta^2 \sigma^2} + \frac {1} {T} \sum _{r=0} ^{R-1} \sum _{k=0} ^{\tau-1} 3 L^2 \Eb {\hat \Xi_{r \tau + k}} \\
     &\refLE{eq:normal distance lemma LocalSGD} \frac {2 \Delta} {\eta T} + \eta \frac {L \sigma^2} {n} + { 9 (\tau-1) L^2 \eta^2 \sigma^2} + \frac {27} {2} { (\tau-1)^2 L^2 \eta^2 \zeta^2 } . 
     \end{aligned}
     \] 
     Finally, \eqref{eq:LocalSGD HS rate} follows from the following assignment 
     \[ 
     \begin{aligned}
     \eta := \min & \Bigg\{ \frac 1 L, \frac {1 - {\rho} / L} {6 \rho (\tau-1)} , \sqrt {\frac {2\Delta n} {L \sigma^2 T}}, \left( \frac {4\Delta} {27 L^2(\tau-1)^2\zeta^2 T} \right) ^{\frac 1 3}, \left( \frac {2\Delta} {9 L^2(\tau-1)\sigma^2 T} \right) ^{\frac 1 3} \Bigg\} .
     \end{aligned}
     \] 
 \end{proof}

 \subsubsection{Proof of Theorem~\ref{thm:LocalSGD faster convex}}
 
 \begin{lemma}[\cf Lemma~7 in~\cite{woodworth2020minibatch}]
     If all $f_i$'s are convex, for $\eta \le \frac {1} {2 L}$, \algname{LocalSGD} ensures
     \begin{equation}
     \label{eq:localsgd descent convex}
     \Eb {f(\bar \xx_t) - f^*} \le \frac 1 \eta \E \norm {\bar \xx_t - \xx^*} ^2 - \frac 1 \eta \E \norm {\bar \xx_{t+1} - \xx^*} ^2 + 3 \frac {\eta \sigma^2} {n} + {2L} \Eb {\Xi_t} . 
     \end{equation}
 \end{lemma}
 
 \begin{lemma}
 \label{lemma:weak convexity distance lemma LocalSGD}
     Under \cref{assumption:GS,assumption:weak convexity}, if all $f_i$'s are convex, for 
     \(
     \eta \le \frac 2 L, 
     \)
     \algname{LocalSGD} ensures 
     \begin{equation}
     \label{eq:weak convexity distance lemma LocalSGD}
     \Eb {\Xi_t} \le 3k \cdot \left[ (\tau-1) \eta^2 \zeta^2 + \eta^2 \sigma^2 \right] ,
     \end{equation}
     where $k = t - r(t)$.
 \end{lemma}
 \begin{proof} 
     Let $\gamma := \frac 1 {\tau-1}$. We prove  
     \begin{equation}
     \label{eq:faster induction distance bound}
     \Eb {\Xi_{t}} \le (1+\gamma) ^{k} \cdot k \cdot \left[ (\tau-1) \eta^2 \zeta^2 + \eta^2 \sigma^2 \right]
     \end{equation}
     by induction on $k$. 
     \begin{enumerate}
     \item \eqref{eq:faster induction distance bound} obviously holds for $k = 0$. 
     \item We assume \eqref{eq:faster induction distance bound} holds for $k$, where $k \le \tau-2$. Then for $\eta \le \frac 2 L$, 
     \begin{equation}
     \label{eq:variance trick in faster}
     \begin{aligned}
         &\quad \Eb {\Xi_{t+1}} \\
         &\refLE{eq:variance relax} \frac 1 n \sum _{i=1} ^n \E \norm { \xx_{t}^i - \eta \gg_t^i - \left( \bar \xx_t - \eta \nabla f(\bar \xx_t) \right) } ^2 \\
         &\le \frac 1 n \sum _{i=1} ^n \E \norm { \xx_{t}^i - \eta \nabla f_i (\xx_t^i) - \left( \bar \xx_t - \eta \nabla f(\bar \xx_t) \right) } ^2 + \eta^2 \sigma^2 \\
         &\refLE{eq:GS bound} \frac 1 n \sum _{i=1} ^n (1+\gamma) \E \norm { \xx_{t}^i - \eta \nabla f_i (\xx_t^i) - \left( \bar \xx_t - \eta \nabla f_i(\bar \xx_t) \right) } ^2 + (1+\gamma^{-1}) \eta^2 \zeta^2 + \eta^2 \sigma^2 \\
         &\refLE{eq:weak convexity contraction} (1+\gamma) \Eb {\Xi_t} + (1+\gamma^{-1}) \eta^2 \zeta^2 + \eta^2\sigma^2 ,
     \end{aligned}
     \end{equation}
     and by induction hypothesis, we have 
     \[
     \begin{aligned}
     &\quad \Eb {\Xi_{t+1}} \\
     &\refLE{eq:faster induction distance bound} \left( 1 + \gamma \right) ^{k+1} \cdot k \left[ (\tau-1) \eta^2 \zeta^2 + \eta^2 \sigma^2 \right] + (1+\gamma^{-1}) \eta^2 \zeta^2 + \eta^2\sigma^2 \\
     &= \left( 1 + \gamma \right) ^{k+1} \cdot k \left[ (\tau-1) \eta^2 \zeta^2 + \eta^2 \sigma^2 \right] + \tau \eta^2 \zeta^2 + \eta^2\sigma^2 \\
     &\le \left( 1 + \gamma \right) ^{k+1} \cdot (k+1) \left[ (\tau-1) \eta^2 \zeta^2 + \eta^2 \sigma^2 \right] .
     \end{aligned}
     \]
     \end{enumerate}
     Now, we have proven \eqref{eq:faster induction distance bound}, and then \eqref{eq:weak convexity distance lemma LocalSGD} follows directly from the fact that 
     \[
     (1+\gamma) ^{k} \le (1+\gamma) ^{\tau-1} < e < 3.
     \]
 \end{proof} 
 
 Now we are ready to prove Theorem~\ref{thm:LocalSGD faster convex}. 
 \begin{proof}[Proof of Theorem~\ref{thm:LocalSGD faster convex}]
     For $\eta \le \frac 1 {2L}$, multiplying both sides of \cref{eq:localsgd descent convex} by $\frac 1 T$ and summing over $t$ from $0$ to $T-1$, we have 
     \[
     \begin{aligned}
     \frac 1 T \sum _{t=0} ^{T-1} \Eb {f(\bar \xx_t)} - f^* &\le \frac {D^2} {\eta T} + 3 \eta \frac {\sigma^2} {n} + \frac {2L} T \sum _{r=0} ^{R-1} \sum _{k=0} ^{\tau-1} \Eb {\Xi_{r\tau +k}} \\
     &\refLE{eq:weak convexity distance lemma LocalSGD} \frac {D^2} {\eta T} + 3 \eta \frac {\sigma^2} {n} + \left. 3 (\tau-1)^2 L \eta^2 \zeta^2 + 3 (\tau-1) L \eta^2 \sigma^2 \right. .
     \end{aligned}
     \]
     Finally, \eqref{eq:LocalSGD convex faster rate} follows from the following assignment: 
     \[
     \eta := \min \left\{ \frac {1} {2L}, \sqrt { \frac {nD^2} {3 \sigma^2 T} }, \left( \frac {D^2} {3L(\tau-1)^2 \zeta^2 T} \right) ^{\frac 1 3}, \left( \frac {D^2} {3L(\tau-1) \sigma^2 T} \right) ^{\frac 1 3} \right\} .
     \]
 \end{proof}
 
 \subsubsection{Proof of Theorem~\ref{thm:LocalSGD improved conditioning}} 
 
 \begin{lemma}
 \label{lemma:localsgd variance lemma}
     For 
     \(
     \eta \le \frac {1} {2 L (\tau-1)} , 
     \)
     \algname{LocalSGD} ensures 
     \begin{equation}
     \label{eq:normal distance lemma LocalSGD distance to noiseless}
     \E \norm { \xx_t^i - \hat \xx_t^i }^2 \le 3 k \eta^2 \sigma^2 , 
     \end{equation}
     where $k = t - r(t)$. 
 \end{lemma}
 \begin{proof}
     Let $\gamma := \frac 1 {2(\tau-1)}$. We prove  
     \begin{equation}
     \label{eq:normal induction distance bound distance to noiseless}
     \E \norm { \xx_t^i - \hat \xx_t^i }^2 \le (1+\gamma) ^{2\cdot k} \cdot k \cdot \eta^2 \sigma^2 
     \end{equation}
     by induction on $k$. 
     \begin{enumerate}
     \item \eqref{eq:normal induction distance bound distance to noiseless} obviously holds for $k = 0$. 
     \item We assume \eqref{eq:normal induction distance bound distance to noiseless} holds for $k$, where $k \le \tau-2$. Then, since $\eta \le \frac 1 {2L (\tau-1)}$, 
     \begin{equation}
     \label{eq:variance trick in conditioning distance to noiseless}
     \begin{aligned}
         \E \norm { \xx_{t+1}^i - \hat \xx_{t+1}^i }^2 
         &\le \E \norm { \xx_{t}^i - \eta \nabla f_i (\xx_t^i) - \hat \xx_{t}^i + \eta \nabla f_i (\hat \xx_t^i) }^2 + \eta^2 \sigma^2 \\
         &\refLE{eq:smooth contraction} \left( 1 + L \eta \right)^2 \E \norm { \xx_t^i - \hat \xx_t^i }^2 + \eta^2 \sigma^2 \\ 
         &\le (1+\gamma)^2 \cdot \E \norm { \xx_t^i - \hat \xx_t^i }^2 + \eta^2 \sigma^2 , 
     \end{aligned}
     \end{equation}
     and by induction hypothesis, we have 
     \[
     \begin{aligned}
     \E \norm { \xx_{t+1}^i - \hat \xx_{t+1}^i }^2 
     &\refLE{eq:normal induction distance bound distance to noiseless} \left( 1 + \gamma \right) ^{2(k+1)} \cdot k \cdot \eta^2 \sigma^2 + \eta^2 \sigma^2 \\
     &\le \left( 1 + \gamma \right) ^{2(k+1)} \cdot (k+1) \eta^2 \sigma^2 .
     \end{aligned}
     \]
     \end{enumerate}
     Now, we have proven \eqref{eq:normal induction distance bound distance to noiseless}, and then \eqref{eq:normal distance lemma LocalSGD distance to noiseless} follows directly from the fact that 
     \[
     (1+\gamma) ^{2 \cdot k} \le (1+\gamma) ^{2(\tau-1)} < e < 3.
     \]
 \end{proof}
 
 \begin{lemma}
 \label{lemma:normal distance lemma LocalSGD noiseless}
     Under \cref{assumption:GS}, for 
     \(
     \eta \le \frac {1} {3 L (\tau-1)} , 
     \)
     \algname{LocalSGD} ensures 
     \begin{equation}
     \label{eq:normal distance lemma LocalSGD noiseless}
     \hat \Xi_t \le 9 k (\tau-1) \eta^2 \zeta^2 , 
     \end{equation}
     where $k = t - r(t)$.
 \end{lemma}
 \begin{proof}
     Let $\gamma := \frac 1 {3(\tau-1)}$. We prove  
     \begin{equation}
     \label{eq:normal induction distance bound noiseless}
     {\hat \Xi_{t}} \le (1+\gamma) ^{3\cdot k} \cdot k \cdot 3 (\tau-1) \eta^2 \zeta^2 
     \end{equation}
     by induction on $k$. 
     \begin{enumerate}
     \item \eqref{eq:normal induction distance bound noiseless} obviously holds for $k = 0$. 
     \item We assume \eqref{eq:normal induction distance bound noiseless} holds for $k$, where $k \le \tau-2$. Then, since $\eta \le \frac 1 {3L (\tau-1)}$, 
     \begin{equation}
     \label{eq:variance trick in conditioning noiseless}
     \begin{aligned}
         {\hat \Xi_{t+1}} 
         &\refLE{eq:variance trick in localsgd} \frac 1 n \sum _{i=1} ^n (1 + \gamma) \norm { \hat \xx_{t}^i - \eta \nabla f_i(\hat \xx_t^i) - \left( \hat \xx_t - \eta \nabla f_i(\hat \xx_t) \right) } ^2 + (1+\gamma^{-1}) \eta^2 \zeta^2 \\
         &\refLE{eq:smooth contraction} (1+\gamma) \cdot \left( 1 + L \eta \right)^2 {\hat \Xi_t} + (1+\gamma^{-1}) \eta^2 \zeta^2 \\ 
         &\le (1+\gamma)^3 \cdot {\hat \Xi_t} + (1+\gamma^{-1}) \eta^2 \zeta^2 . 
     \end{aligned}
     \end{equation}
     and by induction hypothesis, we have 
     \[
     \begin{aligned}
    {\hat \Xi_{t+1}} 
     &\refLE{eq:normal induction distance bound noiseless} \left( 1 + \gamma \right) ^{3(k+1)} \cdot k \cdot 3 (\tau-1) \eta^2 \zeta^2 + (1+\gamma^{-1}) \eta^2 \zeta^2 \\
     &= \left( 1 + \gamma \right) ^{3(k+1)} \cdot k \cdot 3 (\tau-1) \eta^2 \zeta^2 + (3\tau-2) \eta^2 \zeta^2 \\
     &\le \left( 1 + \gamma \right) ^{3(k+1)} \cdot (k+1) \cdot 3 (\tau-1) \eta^2 \zeta^2 .
     \end{aligned}
     \]
     \end{enumerate}
     Now, we have proven \eqref{eq:normal induction distance bound noiseless}, and then \eqref{eq:normal distance lemma LocalSGD noiseless} follows directly from the fact that 
     \[
     (1+\gamma) ^{3 \cdot k} \le (1+\gamma) ^{3(\tau-1)} < e < 3.
     \]
 \end{proof}

 Now we are ready to prove Theorem~\ref{thm:LocalSGD improved conditioning}. 
 \begin{proof}[Proof of Theorem~\ref{thm:LocalSGD improved conditioning}]
     From Lemma~\ref{lemma:refined upper bound of Q}, we have 
     \begin{equation}
     \label{eq:refined Q_t bound in localsgd noiseless}
     \begin{aligned}
     {\hat Q_t} \le 8 \bardelta^2 {\hat \Xi_t} + \frac {\cM^2} 2 {\hat \Xi_t^2} .  
     \end{aligned}
    \end{equation}

     Hence, for $\eta \le \frac 1 {3L (\tau-1)}$, 
     \[
     \begin{aligned}
     &\quad \frac 1 T \sum _{t=0} ^{T-1} \E \norm {\nabla f(\bar \xx_t)}^2 \\
     &\refLE{eq:localsgd descent lemma} \frac {2 \Delta} {\eta T} + \eta \frac {L \sigma^2} {n} + \frac {1} {T} \sum _{r=0} ^{R-1} \sum _{k=0} ^{\tau-1} \left( {6 L^2 \frac 1 n \sum _{i=1} ^n \Eb { \norm { \xx_{r \tau + k}^i - \hat \xx_{r \tau + k}^i }^2 } + 3 \Eb {\hat Q _{r \tau + k}}} \right) \\
     &\stackrel {\eqref{eq:normal distance lemma LocalSGD distance to noiseless} \eqref{eq:refined Q_t bound in localsgd noiseless}} {\le} \frac {2 \Delta} {\eta T} + \eta \frac {L \sigma^2} {n} + { 18 (\tau-1) L^2 \eta^2 \sigma^2} + \frac {1} {T} \sum _{r=0} ^{R-1} \sum _{k=0} ^{\tau-1} \left( 24 \bardelta^2 \Eb {\hat \Xi_{r \tau + k}} + \frac 3 2 \cM^2 \Eb {\hat \Xi_{r \tau + k}^2} \right) \\
     &\refLE{eq:normal distance lemma LocalSGD noiseless} \frac {2 \Delta} {\eta T} + \eta \frac {L \sigma^2} {n} + { 18 (\tau-1) L^2 \eta^2 \sigma^2} + 108 { (\tau-1)^2 \bardelta^2 \eta^2 \zeta^2 } + \frac {81} {4} \cM^2 (\tau-1)^3 (2\tau-1) \eta^4 \zeta^4 . 
     \end{aligned}
     \] 
     Note that in the last inequality, we can upper bound $\Eb {\hat \Xi_{r \tau + k}^2}$ by $(9k(\tau-1)\eta^2\zeta^2)^2$, because the bound in \eqref{eq:normal distance lemma LocalSGD noiseless} is deterministic. (This is why it is necessary to construct the noiseless sequence in our proof. Otherwise, we cannot upper bound $\Eb {\Xi_{r \tau + k}^2}$ likewise.) 
 
     Finally, \eqref{eq:LocalSGD HS rate} follows from the following assignment 
     \[ 
     \begin{aligned}
     \eta := \min & \Bigg\{ \frac 1 L, \frac 1 {3L(\tau-1)}, \sqrt {\frac {2\Delta n} {L \sigma^2 T}}, \left( \frac {\Delta} {54 \bardelta^2 (\tau-1)^2 \zeta^2 T } \right) ^{\frac 1 3}, \left( \frac {\Delta} {9 L^2 (\tau-1) \sigma^2 T } \right) ^{\frac 1 3} , \\
     &\qquad \left( \frac {8 \Delta} {81 \cM^2 (\tau-1)^3 (2\tau-1) \zeta^4 T } \right) ^{\frac 1 5} \Bigg\} .
     \end{aligned}
     \] 
 \end{proof}

 \subsection{The Analysis of SCAFFOLD}
 \label{subsec:the analysis of VR-LocalSGD}
 
 In section~\ref{subsec:the analysis of VR-LocalSGD}, to simplify the notations, for $t=r\tau+k$, $r \in [0, R-1]$, $k \in [0, \tau-1]$, we denote
 \[
 r(t) = r\tau, \quad \xx _{(t)} ^i = \xx _{2r \tau + \tau + k} ^i , \quad \gg _{(t)} ^i = \gg _{2r \tau + \tau + k} ^i ,
 \]
 and
 \[
 \begin{aligned}
 \bar \xx _{(t)} ^i &= \frac 1 n \sum _{i=1} ^n \xx _{(t)} ^i , \\ 
 \Xi _{(t)} &= \frac 1 n \sum _{i=1} ^n \norm { \xx_{(t)}^i - \bar \xx_{(t)}}^2 , \\ 
 Q_{(t)} &= \norm { \frac 1 n \sum _{i=1} ^n \nabla f_i(\xx_{(t)}^i) - \nabla f(\bar \xx_{(t)}) }^2 .
 \end{aligned}
 \] 
 Also, we denote $\bar \xx_{(R\tau)} = \bar \xx_{T}$. 
 
 \begin{lemma}
   For $\eta \le \frac 1 {2L}$, \algname{SCAFFOLD} ensures
   \begin{equation}
   \label{eq:scaffold descent lemma}
   \begin{aligned}
      &\quad \frac 2 T \sum _{t=0}^{R\tau-1} \left[ \E \norm {\nabla f(\bar \xx_{(t)})}^2 + \frac {1} {2} \E \norm { \frac 1 n \sum _{i=1} ^n \nabla f_i(\xx_{(t)}^i) } ^2 \right] \\
      &\le \frac {4 \Delta} {\eta T} + 2\eta \frac {L \sigma^2} {n} + \frac {2} {T} \sum _{r=0} ^{R-1} \sum _{k=0} ^{\tau-1} \Eb {Q_{(r\tau+k)}} . 
   \end{aligned}
   \end{equation}
 \end{lemma}
 \begin{proof}
     For $t=r\tau+k$, $r \in [0, R-1]$, $k \in [0, \tau-1]$, for $\eta \le \frac 1 {2L}$, we have 
     \[
     \begin{aligned}
     &\quad \Eb { f(\bar \xx_{(t+1)}) } \\
     &\le \Eb {f(\bar \xx_{(t)})} - \eta \E \lin { \nabla f(\bar \xx_{(t)}) , \frac 1 n \sum _{i=1} ^n \left[ \gg_{(t)}^i - \hat \gg_{(r(t))}^i + \hat \gg_{(r(t))} \right] } \\
     &\qquad + \eta^2 \frac L 2 \E \norm { \frac 1 n \sum _{i=1} ^n \left[ \gg_{(t)}^i - \hat \gg_{(r(t))}^i + \hat \gg_{(r(t))} \right] }^2 \\
     &= \Eb { f(\bar \xx_{(t)}) } - \eta \E \lin { \nabla f(\bar \xx_{(t)}), \frac 1 n \sum _{i=1} ^n \gg _{(t)}^i } + \frac L 2 \eta^2 \E \norm { \frac 1 n \sum _{i=1} ^n \gg _{(t)}^i }^2 \\
     &\refLE{eq:borel} \Eb { f(\bar \xx_{(t)}) } - \eta \E \lin { \nabla f(\bar \xx_{(t)}) , \frac 1 n \sum _{i=1} ^n \nabla f_i(\xx_{(t)}^i) } + \frac L 2 \eta ^2 \E \norm { \frac 1 n \sum _{i=1} ^n \nabla f_i(\xx_{(t)}^i) } ^2 + \frac {\eta ^2} 2 L \frac {\sigma^2} n \\
     &= \Eb { f(\bar \xx_{(t)}) } - \frac {\eta} 2 \E \norm { \nabla f(\bar \xx_{(t)}) } ^2 - \frac {\eta} 4 \E \norm { \frac 1 n \sum _{i=1} ^n \nabla f_i(\xx_{(t)}^i) } ^2 + \frac {\eta} 2 \Eb { Q_{(t)} } + \frac {\eta^2} 2 L \frac {\sigma^2} n . 
     \end{aligned}
     \] 
     Summing over $k$ from $0$ to $\tau-1$ and over $r$ from $0$ to $R-1$, and multiplying both sides by $\frac 4 {\eta T}$, we have 
     \[ 
     \begin{aligned}
     &\quad \frac 2 T \sum _{t=0}^{R\tau-1} \left[ \E \norm {\nabla f(\bar \xx_{(t)})}^2 + \frac {1} {2} \E \norm { \frac 1 n \sum _{i=1} ^n \nabla f_i(\xx_{(t)}^i) } ^2 \right] \\
     &\le \frac {4} {\eta T} \E \left[ f(\bar \xx_{(0)}) - f(\bar \xx_{(T)}) \right] + 2\eta \frac {L \sigma^2} {n} + \frac {2} {T} \sum _{r=0} ^{R-1} \sum _{k=0} ^{\tau-1} \Eb {Q_{(r\tau+k)}} \\ 
     &\le \frac {4 \Delta} {\eta T} + 2\eta \frac {L \sigma^2} {n} + \frac {2} {T} \sum _{r=0} ^{R-1} \sum _{k=0} ^{\tau-1} \Eb {Q_{(r\tau+k)}} . 
     \end{aligned}
     \]
 \end{proof}
 
 \begin{lemma}
     Under \cref{assumption:mean HS,assumption:weak convexity}, for $\gamma = \frac 1 {3(\tau-1)}$ and $\eta \le \frac {(1-\rho/L)} {2 \rho} \gamma$, we have 
     \begin{equation}
     \label{eq:scaffold recursion}
     \begin{aligned}
     \Eb {\Xi_{(t)}} 
     &\le 3 \gamma^{-2} (1+\gamma^{-1}) \eta^4 \delta^2 \sum _{l=0} ^{k-1} \E \norm { \frac 1 n \sum _{i=1} ^n \nabla f_i (\xx_{(r(t)+l)}^i) } ^2 \\
     &\qquad + 3 \left( \frac { 1 + \gamma^{-1} } { \tau } + 1 \right) \eta^2 k \sigma^2 + 3 \gamma^{-2} \eta^4 k \delta^2 \frac {\sigma^2} n , 
     \end{aligned}
     \end{equation}
     where $k = t-r(t)$. 
 \end{lemma}
 \begin{proof}
 For $\gamma > 0$, when $t+1 > r(t+1)$, in view of 
 \begin{equation}
 \label{eq:Xi_t recursion}
 \begin{aligned}
 &\quad \Eb { \Xi_{(t+1)} } \\
 &\refLE{eq:variance relax} \frac 1 n \sum _{i=1} ^n \E \norm { \xx_{(t+1)}^i - \left( \bar \xx_{(t)} - \eta \nabla f(\bar \xx_{(t)}) \right) } ^2 \\
 &\refLE{eq:scaffold_update} \frac 1 n \sum _{i=1} ^n \E \norm { \xx_{(t)}^i - \eta \left( \gg_{(t)}^i - \hat \gg_{(r(t))}^i + \hat \gg_{(r(t))} \right) - \bar \xx_{(t)} + \eta \nabla f(\bar \xx_{(t)}) } ^2 \\
 &\refLE{eq:borel} \frac 1 n \sum _{i=1} ^n \E \norm { \xx_{(t)}^i - \eta \left( \nabla f_i(\xx_{(t)}^i) - \hat \gg_{(r(t))}^i + \hat \gg_{(r(t))} \right) - \bar \xx_{(t)} + \eta \nabla f(\bar \xx_{(t)}) } ^2 + \eta^2 \sigma^2 \\
 &\le \left( 1 + \gamma \right) \frac 1 n \sum _{i=1} ^n \E \norm { \xx_{(t)}^i - \eta \left( \nabla f_i(\xx_{(t)}^i) - \nabla f_i(\bar \xx_{(r(t))}) + \nabla f(\bar \xx_{(r(t))}) \right) - \bar \xx_{(t)} + \eta \nabla f(\bar \xx_{(t)}) } ^2 \\
 &\qquad + \left( 1 + \gamma ^{-1} \right) \eta^2 \E \norm { \hat \gg_{(r(t))}^i - \hat \gg_{(r(t))} - \nabla f_i(\bar \xx_{(r(t))}) + \nabla f(\bar \xx_{(r(t))}) } ^2 + \eta^2 \sigma^2 \\
 &\refLE{eq:borel} \left( 1 + \gamma \right) \frac 1 n \sum _{i=1} ^n \E \norm { \xx_{(t)}^i - \eta \left( \nabla f_i(\xx_{(t)}^i) - \nabla f_i(\bar \xx_{(r(t))}) + \nabla f(\bar \xx_{(r(t))}) \right) - \bar \xx_{(t)} + \eta \nabla f(\bar \xx_{(t)}) } ^2 \\ 
 &\qquad + \left( \frac { 1 + \gamma^{-1} } { \tau } + 1 \right) \eta^2 \sigma^2 , 
 \end{aligned}
 \end{equation}
 and, for $\eta \le \frac {(1-\rho/L)} {2 \rho} \gamma$, 
 \begin{equation}
 \label{eq:Xi_t variance reduction}
 \begin{aligned}
     &\quad \frac 1 n \sum _{i=1} ^n \E \norm { \xx_{(t)}^i - \eta \left( \nabla f_i(\xx_{(t)}^i) - \nabla f_i(\bar \xx_{(r(t))}) + \nabla f(\bar \xx_{(r(t))}) \right) - \bar \xx_{(t)} + \eta \nabla f(\bar \xx_{(t)}) } ^2 \\
     &\le \frac 1 n \sum _{i=1} ^n \E \bigg[ (1 + \gamma) \norm { \xx_{(t)}^i - \eta \nabla f_i (\xx_{(t)}^i) - \bar \xx_{(t)} + \eta \nabla f_i(\bar \xx_{(t)}) } ^2 \\
     &\qquad + (1 + \gamma ^{-1}) \eta^2 \norm { \nabla f_i (\bar \xx_{(r(t))}) - \nabla f(\bar \xx_{(r(t))}) - \nabla f_i (\bar \xx_{(t)}) + \nabla f (\bar \xx_{(t)}) } ^2 \bigg] \\
     &\refLE{eq:mean HS bound} \frac 1 n \sum _{i=1} ^n \E \bigg[ (1 + \gamma) \norm { \xx_{(t)}^i - \eta \nabla f_i (\xx_{(t)}^i) - \bar \xx_{(t)} + \eta \nabla f_i(\bar \xx_{(t)}) } ^2 + (1 + \gamma ^{-1}) \eta^2 \delta^2 \norm { \bar \xx_{(t)} - \bar \xx_{(r(t))} } ^2 \bigg] \\
     &\refLE{eq:weak convexity contraction} (1 + \gamma) \left( 1 + \frac {2L} {L-\rho} \cdot \rho \eta \right) \Eb {\Xi_{(t)}} + (1+ \gamma ^{-1}) \eta^2 \delta^2 \E \norm { \bar \xx_{(t)} - \bar \xx_{(r(t))} } ^2 \\
     &\le (1+\gamma)^2 \Eb {\Xi_{(t)}} + (1+ \gamma ^{-1}) \eta^2 \delta^2 \E \norm { \bar \xx_{(t)} - \bar \xx_{(r(t))} } ^2 ,  
 \end{aligned}
 \end{equation}
 we have
 \[
 \begin{aligned}
 \label{eq:Xi recursion in SCAFFOLD}
 &\quad \Eb { \Xi_{(t+1)} } \\
 &\le (1+\gamma)^3 \Eb {\Xi_t} + (1+\gamma) (1+\gamma^{-1}) \eta^2 \delta^2 \E \norm { \bar \xx_{(t)} - \bar \xx_{(r(t))} } ^2 + \left( \frac { 1 + \gamma^{-1} } { \tau } + 1 \right) \eta^2 \sigma^2 . 
 \end{aligned}
 \]
 
 Also, when $t+1 > r(t+1)$, 
 \begin{equation}
 \label{eq:averaged move recursion}
 \begin{aligned}
 &\quad \E \norm { \bar \xx_{(t+1)} - \bar \xx_{(r(t+1))} } ^2 \\
 &\refEQ{eq:scaffold_update} \E \norm { \bar \xx_{(t)} - \eta \cdot \frac 1 n \sum _{i=1} ^n \gg_{(t)}^i - \bar \xx_{(r(t))} } ^2 \\
 &\refLE{eq:borel} \E \norm { \bar \xx_{(t)} - \eta \cdot \frac 1 n \sum _{i=1} ^n \nabla f_i (\xx_{(t)}^i) - \bar \xx_{(r(t))} } ^2 + \eta^2 \frac {\sigma^2} n \\
 &\le (1+ \gamma) \E \norm { \bar \xx_{(t)} - \bar \xx_{(r(t))} } ^2 + (1+ \gamma^{-1}) \eta^2 \E \norm { \frac 1 n \sum _{i=1} ^n \nabla f_i (\xx_{(t)}^i) } ^2 + \eta^2 \frac {\sigma^2} n .
 \end{aligned}
 \end{equation}
 
 Now, when $t+1 > r(t+1)$, we have 
 \[
 \begin{aligned}
     &\quad \Eb {\Xi_{(t+1)} + \gamma^{-2} \eta^2 \delta^2 \norm { \bar \xx_{(t+1)} - \bar \xx_{(r(t+1))} } ^2} \\
     &\refLE{eq:Xi_t recursion} (1+\gamma)^3 \cdot \Eb {\Xi_{(t)}} + (1+\gamma)(1+\gamma^{-1}) \eta^2 \delta^2 \E \norm { \bar \xx_{(t)} - \bar \xx_{(r(t))} } ^2 \\ 
     &\qquad + \gamma^{-2} \eta^2 \delta^2 \E \norm { \bar \xx_{(t+1)} - \bar \xx_{(r(t+1))} } ^2 + \left( \frac { 1 + \gamma^{-1} } { \tau } + 1 \right) \eta^2 \sigma^2 \\
     &\refLE{eq:averaged move recursion} (1+\gamma)^3 \cdot \Eb {\Xi_{(t)}} + \left[ (1+\gamma) (1+\gamma^{-1}) + (1+\gamma) \gamma ^{-2} \right] \eta^2 \delta^2 \E \norm { \bar \xx_{(t)} - \bar \xx_{(r(t))} } ^2 \\
     &\qquad + \gamma^{-2} (1+\gamma^{-1}) \eta^4 \delta^2 \E \norm { \frac 1 n \sum _{i=1} ^n \nabla f_i (\xx_{(t)}^i) } ^2 + \left( \frac { 1 + \gamma^{-1} } { \tau } + 1 \right) \eta^2 \sigma^2 + \gamma^{-2} \eta^4 \delta^2 \frac {\sigma^2} n \\
     &\le (1+\gamma)^3 \cdot \Eb {\Xi_{(t)} + \gamma^{-2} \eta^2 \delta^2 \norm { \bar \xx_{(t)} - \bar \xx_{(r(t))} } ^2} \\ 
     &\qquad + \gamma^{-2} (1+\gamma^{-1}) \eta^4 \delta^2 \E \norm { \frac 1 n \sum _{i=1} ^n \nabla f_i (\xx_{(t)}^i) } ^2 + \left( \frac { 1 + \gamma^{-1} } { \tau } + 1 \right) \eta^2 \sigma^2 + \gamma^{-2} \eta^4 \delta^2 \frac {\sigma^2} n . 
 \end{aligned}
 \]
 
     Finally, for $\gamma \le \frac 1 {3(\tau-1)}$, 
     \[
     \begin{aligned}
     &\quad \Eb {\Xi_{(t)} + \gamma^{-2} \eta^2 \delta^2 \norm { \bar \xx_{(t)} - \bar \xx_{(r(t))} } ^2} \\
     &\le (1+\gamma)^3 \cdot \Eb {\Xi_{(t-1)} + \gamma^{-2} \eta^2 \delta^2 \norm { \bar \xx_{(t-1)} - \bar \xx_{(r(t-1))} } ^2} \\ 
     &\qquad + \gamma^{-2} (1+\gamma^{-1}) \eta^4 \delta^2 \E \norm { \frac 1 n \sum _{i=1} ^n \nabla f_i (\xx_{(t-1)}^i) } ^2 + \left( \frac { 1 + \gamma^{-1} } { \tau } + 1 \right) \eta^2 \sigma^2 + \gamma^{-2} \eta^4 \delta^2 \frac {\sigma^2} n \\ 
     &\le \cdots \\
     &\le 0 + e \gamma^{-2} (1+\gamma^{-1}) \eta^4 \delta^2 \sum _{l=0} ^{k-1} \E \norm { \frac 1 n \sum _{i=1} ^n \nabla f_i (\xx_{(r(t)+l)}^i) } ^2 \\
     &\qquad + e \left( \frac { 1 + \gamma^{-1} } { \tau } + 1 \right) \eta^2 k \sigma^2 + e \gamma^{-2} \eta^4 k \delta^2 \frac {\sigma^2} n , 
     \end{aligned}
     \]
     where $k=t-r(t)$. 
 \end{proof}
 
 \subsubsection{Proof of Theorem~\ref{thm:SCAFFOLD speedup}}
 
 \begin{proof}[Proof of Theorem~\ref{thm:SCAFFOLD speedup}]
 For $\eta \le \min \left\{ \frac 1 {2L} , \frac {1-\rho/L} {6 \rho (\tau-1)} \right\}$, 
 \[
 \begin{aligned}
 &\quad \frac 2 T \sum _{t=0}^{R\tau-1} \left[ \E \norm {\nabla f(\bar \xx_{(t)})}^2 + \frac {1} {2} \E \norm { \frac 1 n \sum _{i=1} ^n \nabla f_i(\xx_{(t)}^i) } ^2 \right] \\
 &\refLE{eq:scaffold descent lemma} \frac {4 \Delta} {\eta T} + 2\eta \frac {L \sigma^2} {n} + \frac {2} {T} \sum _{r=0} ^{R-1} \sum _{k=0} ^{\tau-1} \Eb {Q_{(r\tau+k)}} \\
 &\le \frac {4 \Delta} {\eta T} + 2\eta \frac {L \sigma^2} {n} + \frac {2L^2} {T} \sum _{r=0} ^{R-1} \sum _{k=0} ^{\tau-1} \Eb {\Xi_{(r\tau+k)}} \\
 &\refLE{eq:scaffold recursion} \frac {4 \Delta} {\eta T} + 2\eta \frac {L \sigma^2} {n} + \frac {162} {T} \sum _{r=0} ^{R-1} \sum _{k=0} ^{\tau-1} \eta^4 L^2 \delta^2 (\tau-1)^3 \tau \E \norm { \frac 1 n \sum _{i=1} ^n \nabla f_i(\xx_{(r\tau+k)}^i) } ^2 \\
 &\qquad + {6} \eta^2 L^2 (\tau-1) \sigma^2 + \frac {27} {2} \eta^4 L^2 \delta^2 (\tau-1)^3 \frac {\sigma^2} {n} . 
 \end{aligned}
 \]
 Then, for $\eta \le \frac 1 {4 \sqrt {L\delta} \tau}$, 
 \[
 \frac 2 T \sum _{t=0}^{R\tau-1} \E \norm {\nabla f(\bar \xx_{(t)})}^2 \le \frac {4 \Delta} {\eta T} + 2\eta \frac {L \sigma^2} {n} + {6} \eta^2 L^2 (\tau-1) \sigma^2 + \frac {27} {2} \eta^4 L^2 \delta^2 (\tau-1)^3 \frac {\sigma^2} {n} . 
 \]
 Finally, \eqref{eq:improved general SCAFFOLD rate} follows from the following assignment
 \[
 \eta := \min \left\{ \frac 1 {2L} , \frac {1-\rho/L} {6 \rho (\tau-1)}, \frac 1 {4 \sqrt {L\delta} \tau}, \sqrt {\frac {\Delta n} {2L \sigma^2 T}}, \left( \frac {2\Delta} {3L^2(\tau-1)\sigma^2 T} \right) ^{\frac 1 3}  \right\} .
 \]
 \end{proof} 
 
 \subsubsection{Proof of Theorem~\ref{thm:Lipschitz Hessian SCAFFOLD}} 
 
 \begin{proof}[Proof of Theorem~\ref{thm:Lipschitz Hessian SCAFFOLD}]
 For $\eta \le \min \left\{ \frac 1 {2L} , \frac {1-\rho/L} {6 \rho (\tau-1)} \right\}$, 
 \[
 \begin{aligned}
 &\quad \frac 2 T \sum _{t=0}^{R\tau-1} \left[ \E \norm {\nabla f(\bar \xx_{(t)})}^2 + \frac {1} {2} \E \norm { \frac 1 n \sum _{i=1} ^n \nabla f_i(\xx_{(t)}^i) } ^2 \right] \\
 &\refLE{eq:scaffold descent lemma} \frac {4 \Delta} {\eta T} + 2\eta \frac {L \sigma^2} {n} + \frac {2} {T} \sum _{r=0} ^{R-1} \sum _{k=0} ^{\tau-1} \Eb {Q_{(r\tau+k)}} \\
 &\refLE{eq:refined upper bound of Q} \frac {4 \Delta} {\eta T} + 2\eta \frac {L \sigma^2} {n} + \frac {16 \bardelta^2} {T} \sum _{r=0} ^{R-1} \sum _{k=0} ^{\tau-1} \Eb {\Xi_{(r\tau+k)}} + 0 \\
 &\refLE{eq:scaffold recursion} \frac {4 \Delta} {\eta T} + 2\eta \frac {L \sigma^2} {n} + \frac {1296} {T} \sum _{r=0} ^{R-1} \sum _{k=0} ^{\tau-1} \eta^4 \bardelta^2 \delta^2 (\tau-1)^3 \tau \E \norm { \frac 1 n \sum _{i=1} ^n \nabla f_i(\xx_{(r\tau+k)}^i) } ^2 \\
 &\qquad + {48} \eta^2 \bardelta^2 (\tau-1) \sigma^2 + {108} \eta^4 \bardelta^2 \delta^2 (\tau-1)^3 \frac {\sigma^2} {n} . 
 \end{aligned}
 \]
 Then, for $\eta \le \frac 1 {6 \sqrt {\bardelta \delta} \tau}$, 
 \[
 \frac 2 T \sum _{t=0}^{R\tau-1} \E \norm {\nabla f(\bar \xx_{(t)})}^2 \le \frac {4 \Delta} {\eta T} + 2\eta \frac {L \sigma^2} {n} + {48} \eta^2 \bardelta^2 (\tau-1) \sigma^2 + {108} \eta^4 \bardelta^2 \delta^2 (\tau-1)^3 \frac {\sigma^2} {n} . 
 \]
 Finally, \eqref{eq:quadratic SCAFFOLD rate} follows from the following assignment 
 \[
 \eta := \min \left\{ \frac 1 {2L} , \frac {1-\rho/L} {6 \rho (\tau-1)}, \frac 1 {6 \sqrt{ \bardelta \delta} \tau}, \sqrt { \frac {2\Delta n} {L\sigma^2 T} }, \left( \frac {\Delta} {12 \bardelta^2 (\tau-1) \sigma^2 T} \right) ^{\frac 1 3}  \right\} . 
 \]
 \end{proof}

\section{Details and Further Experiments}
\label{sec:expr details}

\subsection{Details of the Experiments}

We use a smoothed variant of Huber loss $h(\cdot)$ in all our experiments. We recall its definition:
\[
h (u) = \begin{cases} 
    \frac 1 2 u^2, & |u| \le 1, \\
    - \frac 1 6 (|u|-1)^3 + \frac 1 2 u^2, & 1 < |u| \le 2, \\
    \frac 3 2 |u| - \frac 7 6, & |u| > 2. 
\end{cases}
\]
In the original Huber loss function, the quadratic loss is used for small errors and linear loss for large errors. While in the smoothed variant, a cubic term is added in the middle range to provide a smooth transition between the quadratic and linear sections, thus ensuring the Lipschitz continuity of the Hessian.

The data is generated as follows. We first generate a diagonal matrix $\oldLambda$ with values ranging from $0$ to a maximal value $L$. Next, we generate a random orthonormal matrix $\mQ$, and compute $\mA = \sqrt {\oldLambda} \mQ \T$. With $\mH = \mA \T \mA$, the quadratic function $\frac 1 2 \xx\T \mH \xx$ has $L$-Lipschitz continuous gradient. Then, we assign the function for each worker: $f_i(x) = \sum _{j=1} ^d h \left( \mA_i(j) \cdot (\xx - \xx^*_i) \right)$, where $\mA_i$ introduces random noise to $\mA$ and $\xx^*_i$ represents a unique optimal point for each worker. The function $h(u)$ behaves quadratically when $|u|$ is small, specifically $f_i(x) = \sum _{j=1} ^d h \left( \mA_i(j) \cdot (\xx - \xx^*_i) \right) \approx \frac 1 2 (\xx - \xx^*_i) \T \mA \T \mA (\xx - \xx^*_i)$. This allows us to approximately control the behavior of the functions. We control the Hessian similarity $\delta$ by adjusting the noise in $\mA_i$, and the gradient similarity $\zeta$ by the variation in $\xx^*_i$. 

We run all the experiments on {Intel Xeon CPU E5-2680v4}. We use three different random seeds $\{111,222,333\}$ for all experiments and report the average results across the runs. 

Codes are available here: \url{https://github.com/riekenluo/Pub_Code_LocalSGD_and_SCAFFOLD}. 

\subsection{Further Experiments}

We run additional experiments to validate the impact of the other parameters on the convergence rates. Following our previous experiment setup, we set $L \approx 1.0$, $\Delta \approx 1.0$, $\sigma \approx 0.01$, and $\tau=50$. To highlight the optimization term, we consider the $10$-th communication round; to highlight the heterogeneity term, we consider the $50$-th communication round. 

The level of non-convexity is controlled by $\lambda$ in our experiments. The regularizer $\frac {x^2} {1+x^2}$ exhibits non-convexity only when $|x| > \frac {\sqrt 3} 3$. Therefore, to validate the dependence on different levels of non-convexity, we need to increase the scales of $\xx^*_i$ in these experiments. 

\begin{figure}[ht]
    \centering
    \subfigure{
        \includegraphics[width=0.27\textwidth]{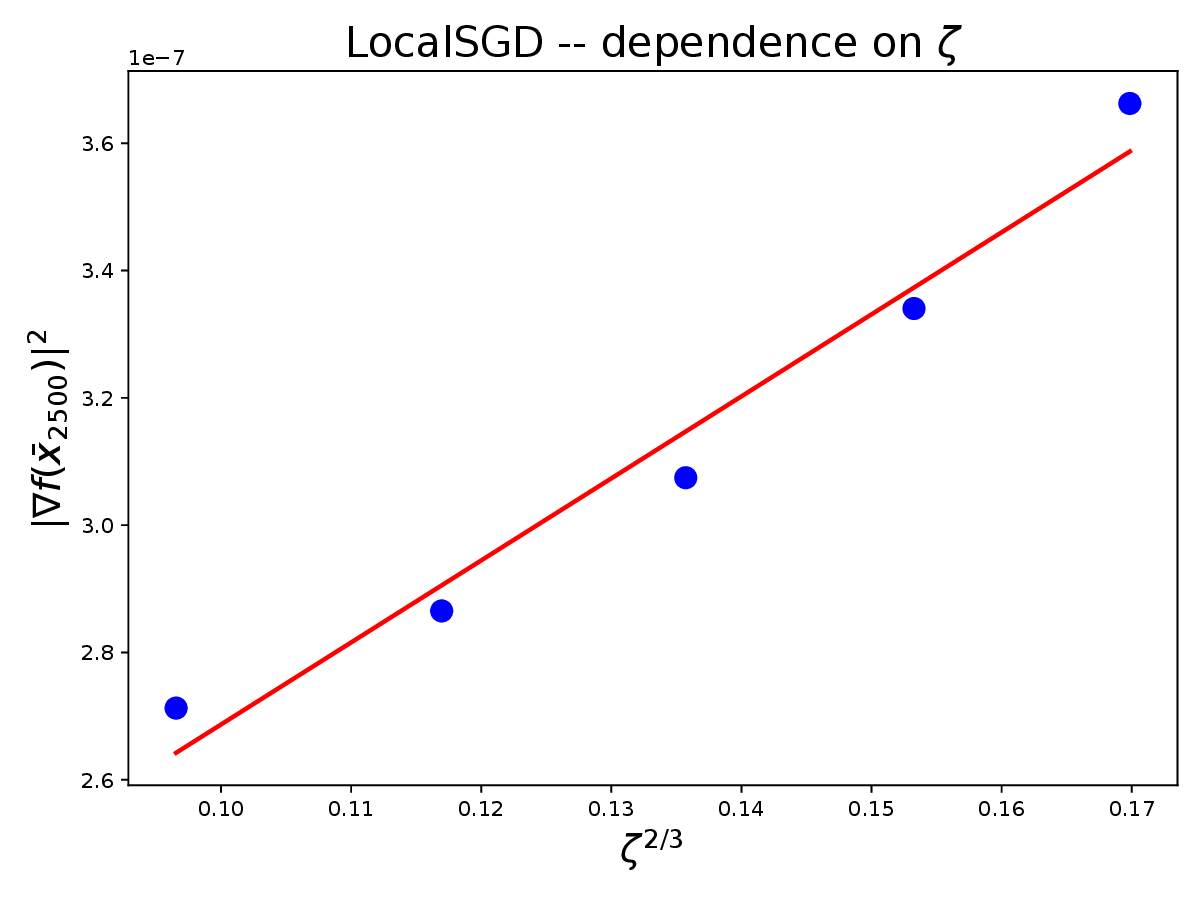} 
    }
    \subfigure{
        \includegraphics[width=0.27\textwidth]{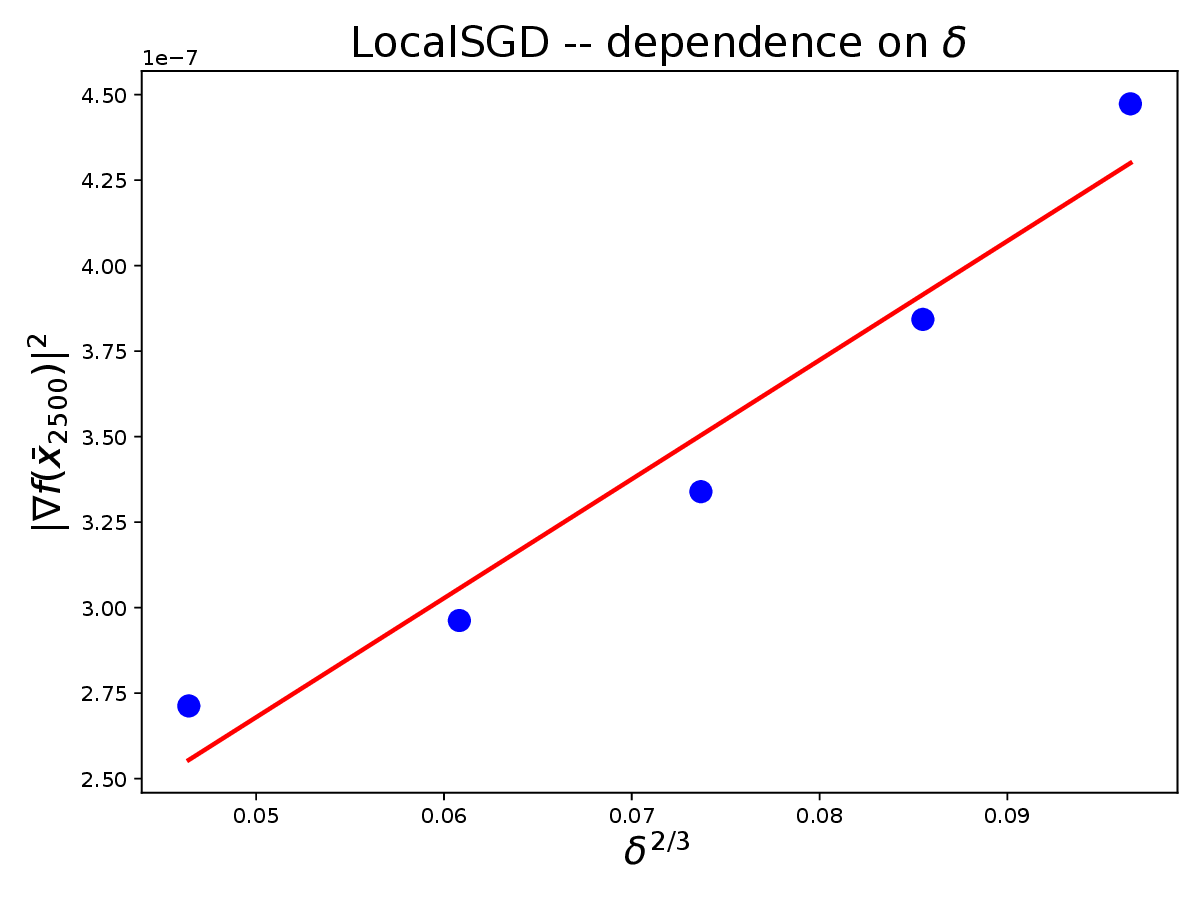}  
    }
    \subfigure{
        \includegraphics[width=0.27\textwidth]{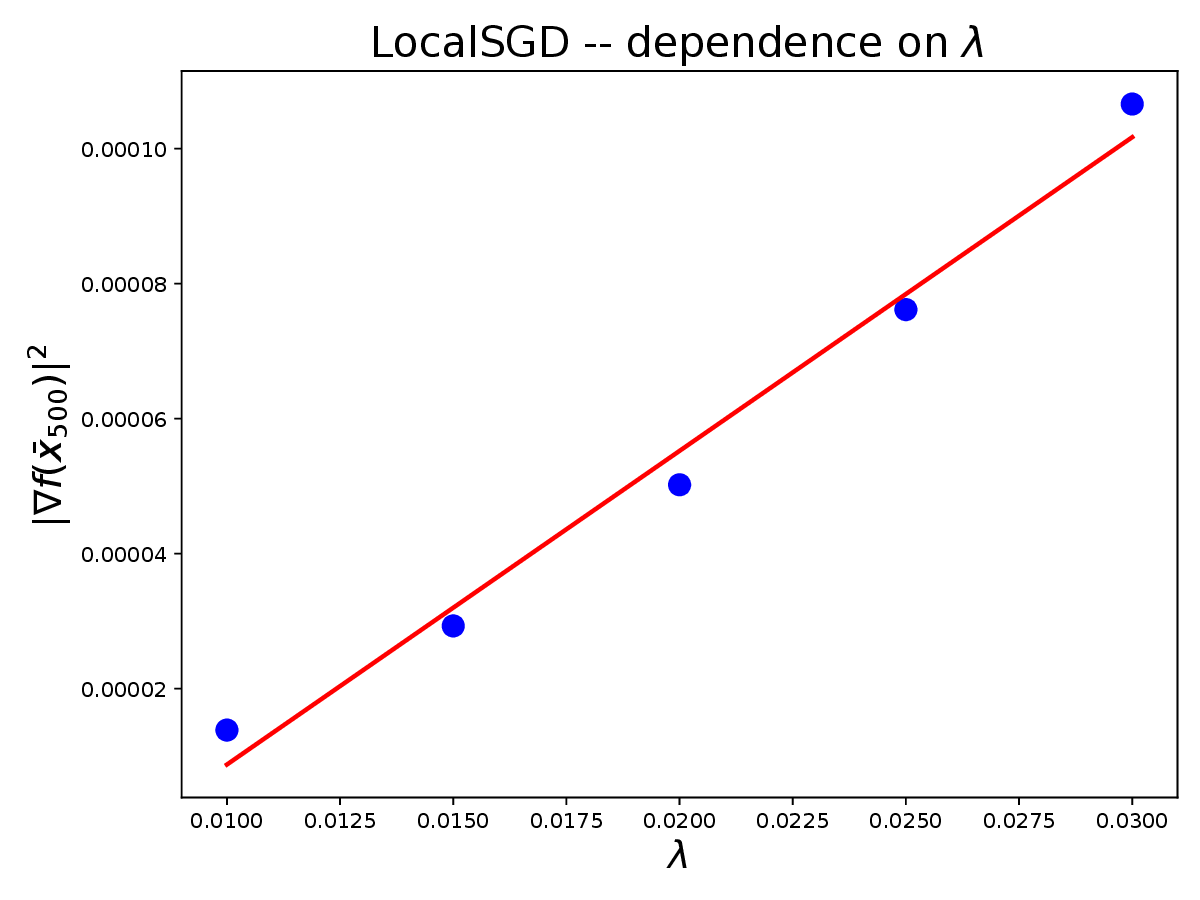}  
    }
    \vspace{-1em}
    \caption{Convergence rates of \algname{LocalSGD} with changing parameters.}
    \label{fig:localsgd dependence}
\end{figure}

For \algname{LocalSGD}, in addition to the experiments in~\cref{sec:expr}, we perform further tests by varying the parameters of gradient similarity, Hessian similarity, and weak convexity. Specifically, from the basic setting $(\zeta,\delta,\lambda)=(0.03,0.01,0.01)$, we change $\zeta$ over $\{0.03, 0.04, 0.05, 0.06, 0.07\}$, $\delta$ over $\{0.01, 0.015, 0.02, 0.025, 0.03\}$, and $\lambda$ over $\{0.01, 0.015, 0.02, 0.025, 0.03\}$. The results of these experiments are plotted in~\cref{fig:localsgd dependence}. 

\begin{figure}[ht]
    \centering
    \subfigure{
        \includegraphics[width=0.27\textwidth]{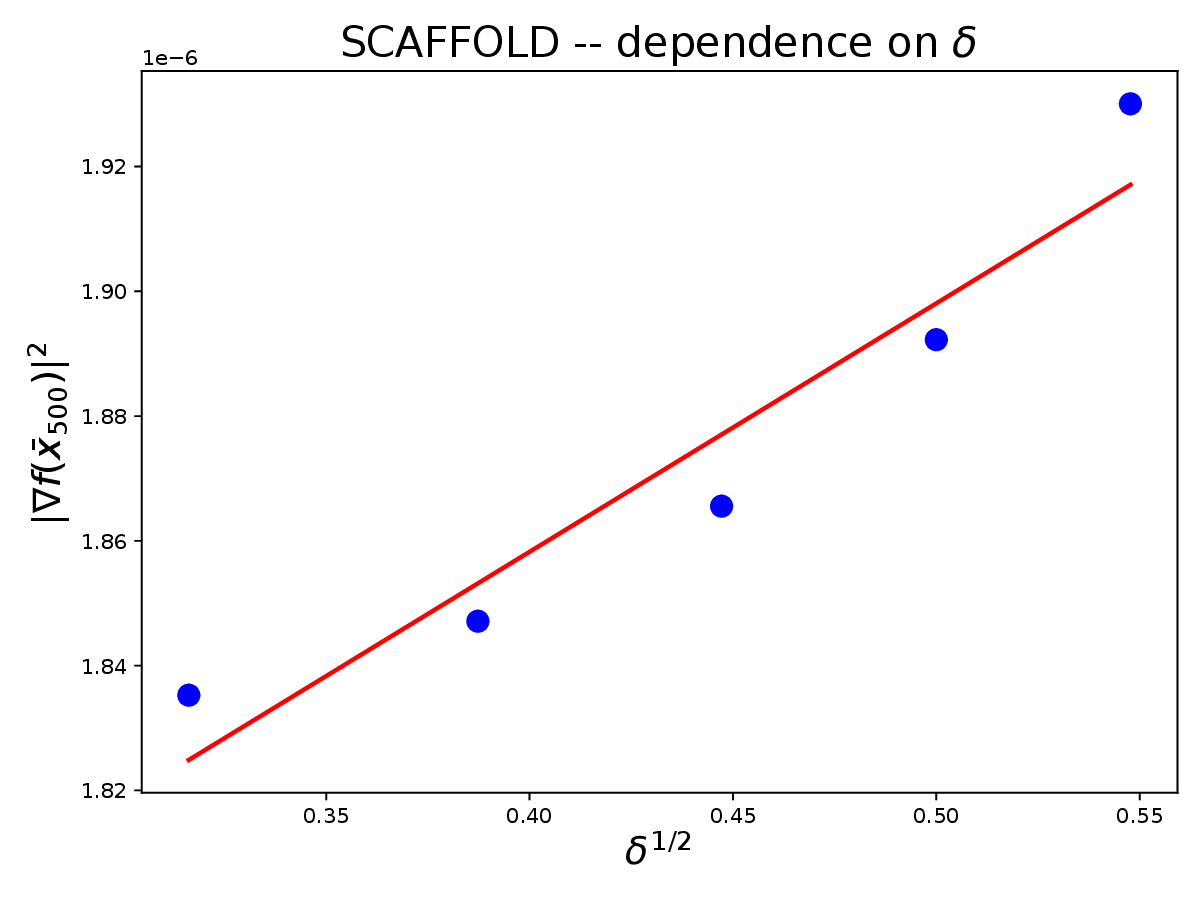}  
    }
    \subfigure{
        \includegraphics[width=0.27\textwidth]{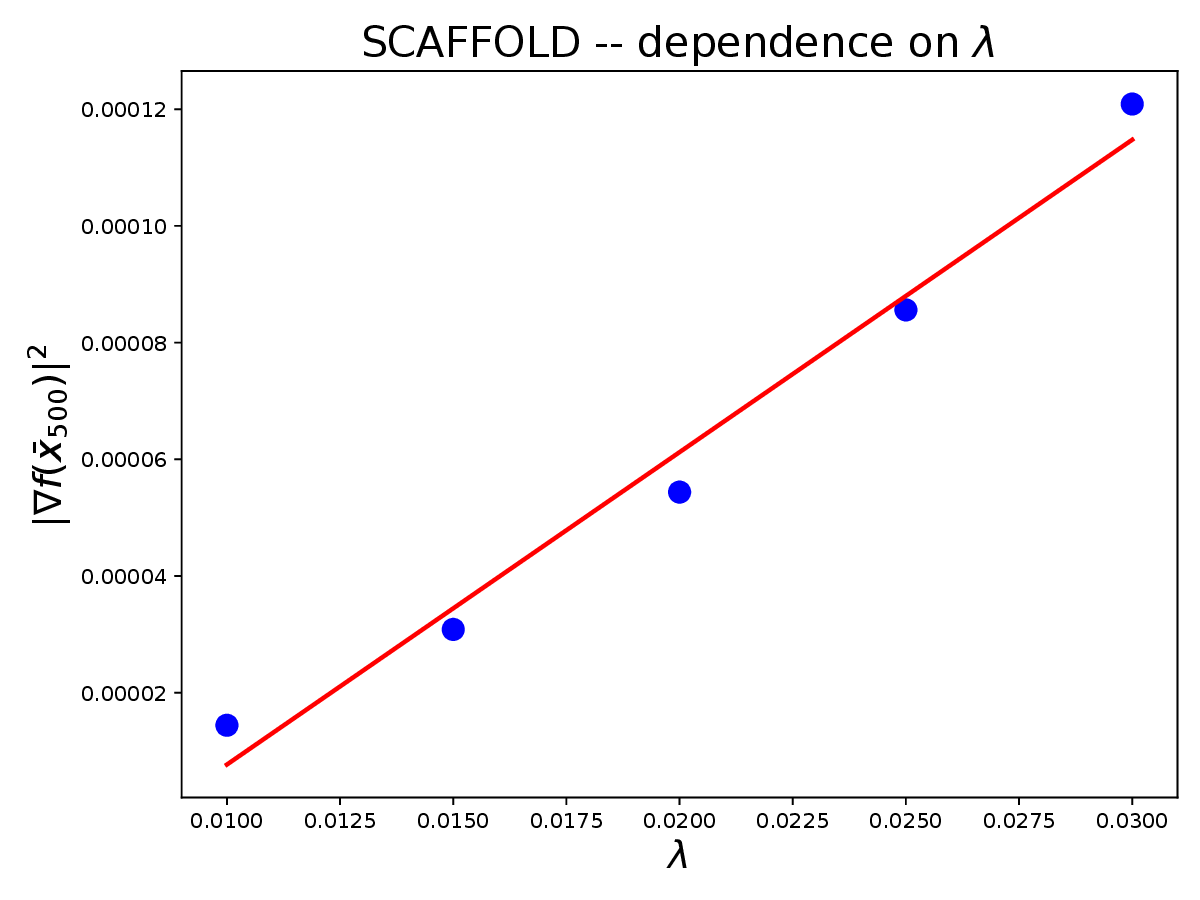}  
    }
    \caption{Convergence rates of \algname{SCAFFOLD} with changing parameters.}
    \label{fig:scaffold dependence}
\end{figure}

For \algname{SCAFFOLD}, in addition to the experiments in~\cref{sec:expr}, we perform further tests by varying the parameters of Hessian similarity and weak convexity. Specifically, from the basic setting $(\delta,\lambda)=(0.1,0.01)$, we change $\delta$ over $\{0.1, 0.15, 0.2, 0.25, 0.3\}$ and $\lambda$ over $\{0.01, 0.015, 0.02, 0.025, 0.03\}$. The results of these experiments are plotted in~\cref{fig:scaffold dependence}. 

The blue dots represent the convergence rates corresponding to each parameter value, and the red lines show the linear fit of these data points. From the figures, it is evident that the convergence rates are positively correlated with these parameters. All these experiments are qualitative in nature, as it is not possible to completely isolate the effects of other additive terms on the convergence rates. 

\section{Discussions on the Related Work on Generalization}

\label{appendix_sec:generalization}

\cite{sefidgaran2024lessons} studies \algname{LocalSGD} when $T$, the number of queries to $\mathcal{SO}$, is fixed. Under the assumption that the loss functions are $\Sigma$-subgaussian, they show that the generalization gap is bounded by $\cO (\Sigma/\sqrt{\tau})$, which decreases as $\tau$ increases. Therefore, for this intermittent communication model, combining the result of~\cite{sefidgaran2024lessons} and the optimization bounds (in this paper or in other related works) yields a preliminary insight for minimizing the true risk, \ie\ to choose the communication interval $\tau$ balancing the empirical risk and the generalization gap. Please refer to~\citet[section~6]{sefidgaran2024lessons} for more discussions.  

\end{document}